\documentstyle[12pt]{article}

\newcommand{\al}{\alpha}

\newcommand{\del}{\delta}
\newcommand{\Del}{\Delta}
\newcommand{\lm}{\lambda}
\newcommand{\th}{\theta}
\newcommand{\ome}{\omega}
\newcommand{\Ome}{\phi}
\newcommand{\1}{{\bf 1}}
\newcommand{\id}{\mbox{\rm id}}
\newcommand{\Hom}{\mbox{\rm Hom\,}}
\newcommand{\End}{\mbox{\rm End\,}}
\newcommand{\Aut}{\mbox{\rm Aut\,}}
\newcommand{\wt}{\mbox{\rm wt}}
\newcommand{\Res}{\mbox{\rm Res\,}}
\newcommand{\Sp}{\mbox{\rm Span\,}}

\newcommand{\ots}{\otimes}
\newcommand{\ts}{\times}
\newcommand{\ops}{\oplus}
\newcommand{\der}{\partial}
\newcommand{\Ops}{\bigoplus}

\newfont{\frak}{eufm10 scaled \magstep1}
\def\h{\mbox{\frak h}}

\newfont{\sfrak}{eufm10}
\def\sh{\mbox{\sfrak h}}

\newfont{\Bbb}{msbm10 scaled \magstep1}
\def\Z{\mbox{\Bbb Z}}
\def\C{\mbox{\Bbb C}}
\def\Q{\mbox{\Bbb Q}}
\def\N{\mbox{\Bbb N}}

\newfont{\sBbb}{msbm10}
\def\sZ{\mbox{\sBbb Z}}
\def\sC{\mbox{\sBbb C}}
\def\sQ{\mbox{\sBbb Q}}
\def\sN{\mbox{\sBbb N}}

\newcommand{\eq}{\begin{eqnarray}}
\newcommand{\eeq}{\end{eqnarray}}
\newcommand{\eqa}{\begin{eqnarray}}
\newcommand{\eeqa}{\end{eqnarray}}
\newcommand{\eqn}{\begin{eqnarray*}}
\newcommand{\eeqn}{\end{eqnarray*}}
\newcommand{\nn}{\nonumber  }

\newtheorem{dfn}{Definition}[section]
\newtheorem{pro}[dfn]{Proposition}
\newtheorem{thm}[dfn]{Theorem}
\newtheorem{theorem}{Theorem}
\newtheorem{lem}[dfn]{Lemma}
\newtheorem{cor}[dfn]{Corollary}

\newtheorem{rem}[dfn]{Remark}

\makeatletter
\@addtoreset{equation}{section}

\makeatother

\def\bl{\begin{lem}\label}
\def\el{\end{lem}}
\def\bt{\begin{thm}\label}
\def\et{\end{thm}}
\def\bp{\begin{pro}\label}
\def\ep{\end{pro}}
\def\br{\begin{rem}\label}
\def\er{\end{rem}}
\def\bc{\begin{cor}\label}
\def\ec{\end{cor}}
\def\bd{\begin{dfn}\label} 
\def\ed{\end{dfn}}
\def\proof{{\it Proof. }}

\renewcommand{\ss}[1]{\mbox{\scriptsize{${#1}$}}}
\def\B{\langle}
\def\K{\rangle}
\def\qed{\hspace*{\fill}{\mbox{$\Box$}}}
\def\swt{\mbox{$\ss\wt$}}
\def\vs{\vspace{3mm}}
\def\M{\mbox{$M(1)$}}
\def\Mp{\mbox{$M(1)^{+}$}}
\def\Mm{\mbox{$M(1)^{-}$}}
\def\Ml{\mbox{$M(1,\lambda)$}}
\def\Mu{\mbox{$M(1,\mu)$}}
\def\Mn{\mbox{$M(1,\nu)$}}
\def\Mt{\mbox{$M(1)(\theta)$}}
\def\Mtp{\mbox{$M(1)(\theta)^{+}$}}
\def\Mtm{\mbox{$M(1)(\theta)^{-}$}}
\def\hh{\mbox{$\hat{\h}$}}
\def\hp{\mbox{$\hat{\h}^{+}$}}
\def\shp{\mbox{$\hat{\sh}^{+}$}}
\def\hm{\mbox{$\hat{\h}[-1]$}}
\def\hmp{\mbox{$\hat{\h}[-1]^{+}$}}
\def\shmp{\mbox{$\hat{\sh}\ss{[-1]^{+}}$}}
\newcommand{\NO}{\,{\raise0.25em\hbox
{$\mathop{\hphantom {\cdot}}
\limits^{_{\circ}}_{^{\circ}}$}}\,}
\def\ch{\mbox{$\rm ch$}}
\newcommand{\fusion}[3]{\mbox{${\ss{\mbox{$\left(\begin{array}
{cr}#3\\ #1\ \ \ #2 \end{array}\right)$}}}$}}
 
\begin{document}

\title{Fusion rules for the free bosonic orbifold \\vertex
operator algebra}
\author{Toshiyuki Abe}
\date{\small {\it Department of Mathematics, Graduate School of Science,\\
Osaka University, Toyonaka, Osaka 560-0043, Japan}\\
{\rm E-mail: sm3002at@ex.ecip.osaka-u.ac.jp}}
\maketitle

\begin{abstract}
Fusion rules among irreducible modules for the free bosonic orbifold
vertex operator algebra are completely determined.
\end{abstract}

\section{Introduction}

We determine the fusion rules for the free bosonic orbifold vertex
operator algebra $\Mp$ which is the fixed point set of the free
bosonic vertex operator algebra $\M$ under an automorphism $\th$ 
of order 2.  

In \cite{Z}, Zhu introduced an associative algebra $A(V)$ (called
Zhu's algebra) associated to a vertex operator algebra $V$. Zhu's
algebra $A(V)$ inherits a part of vertex operator algebra structure
of $V$ which affords many information for $V$-modules. For
example there exists a one to one correspondence between the set of
equivalence classes of irreducible $\N$-gradable $V$-modules and
the set of equivalence classes of irreducible $A(V)$-modules. Later,
in \cite{FZ}, the notion of Zhu's algebra is generalized to an
$A(V)$-bimodule $A(M)$ for an $\N$-gradable $V$-module $M$,
and the fusion rules of rational vertex operator algebras are
completely characterized in terms of these bimodules (see also
\cite{L}). More precisely, for irreducible $V$-modules $M^{i}$
($i=1,2,3)$ with nontrivial top level $M^{i}_{0}$, there exists a
natural injection from
$I$\fusion{M^{1}}{M^{2}}{M^{3}}, which is the space of intertwining
operators of type \fusion{M^{1}}{M^{2}}{M^{3}}, into the dual space
of the contraction
$(M_{0}^{3})^{*}\ots_{A(V)} A(M^{1})\ots_{A(V)} M_{0}^{2}$.
Moreover, if
$V$ is rational, then this map becomes an isomorphism. By using
this isomorphism, fusion rules were calculated for vertex operators
algebra associated with finite dimensional simple Lie algebras
(\cite{FZ}), and for the minimal series (\cite{W}), etc. 

Let $\h$ be a $d$-dimensional complex vector space with a
nondegenerate symmetric bilinear form and let $\hh$ be its
affinization. Then the Fock space
$\M=S(\h\ots t^{-1}\C[t^{-1}])$ is a simple vertex operator algebra
with central charge $d$, and has the automorphism
$\th$ of order $2$ lifted from the map $\h\to\h,$ $h\mapsto -h$ 
(\cite{FLM}). The fixed point set $\Mp$ of $\M$ is also simple vertex
operator algebra and the $-1$-eigenspace $\Mm$ is irreducible
$\Mp$-module.  It is well known that for every $\lm\in\h$, the
Fock space $\Ml=S(\h\ots t^{-1}\C[t^{-1}])\ots \C e^{\lm}$  is
irreducible $\M$-module and the set $\{ \Ml|\lm\in\h\}$ gives all
inequivalent irreducible $\M$-modules. Moreover if $\lm\neq0$,
$\Ml$ is an irreducible $\Mp$-module, and $\Ml$ and $M(1,-\lm)$
are isomorphic each other as $\Mp$-module (see \cite{DM}). In
addition, the $\th$-twisted $M$-module $\Mt$ defined as an induced
module of twisted affine Lie algebra $\hm$ is also $\Mp$-module,
and the $\pm 1$-eigenspaces $\Mt^{\pm}$ for $\th$ give
inequivalent irreducible $\Mp$-modules. It is known that
these irreducible modules $\M^{\pm}$, $\Ml(\simeq M(1,-\lm))$
($\lm\neq0$) and $\Mt^{\pm}$ give all
inequivalent irreducible
$\Mp$-modules  (see \cite{DN} and \cite{DN2}).

In this paper, we only consider the rank one, i.e., $d=1$ case, and
determine  the fusion rules among any triples of irreducible
$\Mp$-modules. One of the main
results in this paper is the following theorem.
\renewcommand{\thetheorem}{}
\begin{theorem}\nn Let $M$, $N$ and $L$ be irreducible
$\Mp$-modules.

\noindent
(1) If $M=\Mp$, then $N_{\ss{\Mp}{ N}}^{L}=\delta_{N,L}$.

\noindent
(2) If $M=\Mm$, then $N_{\ss{\Mm}{ N}}^{L}$ is $0$ or $1$, and
$N_{\ss{\Mm}{ N}}^{L}=1$ if and only if the pair $(N,L)$ is one of
following pairs:
\eqa(\M^{\pm},\M^{\mp}),\ (\Mt^{\pm},\Mt^{\mp}),\ (\Ml,\Mu)\
(\lm^{2}=\mu^{2}).\nn
\eeqa 

\noindent
(3) If $M=\Ml$ ($\lm\neq0$), then $N_{\ss{\Ml}{ N}}^{L}$ is $0$ or $1$,  and
$N_{\ss{\Ml }{N}}^{L}=1$ if and only if the pair $(N,L)$ is one of
following pairs:
\eqa&&(\M^{\pm},\Mu)\ (\lm^{2}=\mu^{2}),\ (\Mu,\Mn)\
(\nu^{2}=(\lm\pm\mu)^{2}),\nn\\
&&(\Mt^{\pm},\Mt^{\pm}),\
(\Mt^{\pm},\Mt^{\mp}).\nn
\eeqa

\noindent
(4) If $M=\Mtp$, then $N_{\ss{\Mtp}{ N}}^{L}$ is $0$ or $1$, and
$N_{\ss{\Mtp}{ N}}^{L}=1$ if and only if the pair $(N,L)$ is one of
following pairs:
\eqa&&(\M^{\pm},\Mt^{\pm}),\ (\Ml,\Mt^{\pm}).\nn
\eeqa

\noindent
(5) If $M=\Mtm$, then $N_{\ss{\Mtm} N}^{L}$ is $0$ or $1$, and
$N_{\ss{\Mtm} N}^{L}=1$ if and only if the pair $(N,L)$ is one of
following pairs:
\eqa&&(\M^{\pm},\Mt^{\mp}),\ (\Ml,\Mt^{\pm}).\nn
\eeqa
\end{theorem}

Let $M^{i}$ ($i=1,2,3$) be irreducible $\Mp$-modules. Then the
classification result of irreducible $\Mp$-modules in \cite{DN} and
formal characters for irreducible $\Mp$-modules show that the
fusion rule
$N_{M^{1}M^{2}}^{M^{3}} $ is invariant under the any permutations
of $\{1,2,3\}$. In more detail, the explicit forms of formal characters
 of irreducible $\Mp$-modules implies that two irreducible $\Mp$-modules
with same formal characters are isomorphic each other, and in
particular that every irreducible $\Mp$-module is isomorphic to 
its contragredient module. Then the above symmetry of fusion rules
 follows from the fact that $N_{MN}^{L}=N_{NM}^{L}=N_{ML'}^{N'}$
hold for modules $M$, $N$, $L$ of a vertex operator algebra, where
$N'$ and $L'$ are contragredient modules of $N$ and $L$
respectively. Next, we prove that
$A(M^{1})$ is generated as $A(\Mp)$-bimodule by at most two
elements which are images of singular vectors of $M^{1}$ viewed as
a module for Virasoro algebra. This is obtained by using the fact
that the $\Mp$ is generated by the Virasoro element and
a singular vector of weight $4$ (see \cite{DG}). Further using
Frenkel-Zhu injection, we prove that the fusion rule
$N_{M^{1}M^{2}}^{M^{3}}$ is less than 2. More detailed study of
the contraction $(M_{0}^{3})^{*}\ots_{A(V)} A(M^{1})\ots_{A(V)}
M_{0}^{2}$ of $A(\Mp)$-modules implies that the fusion rule
$N_{M^{1}M^{2}}^{M^{3}}$ is in fact less than $1$. In \cite{FLM} and
\cite{DN}, the nontrivial intertwining operator of type
\fusion{\Ml}{\Mu}{M(1,\lm+\mu)} was constructed for every
$\lm,\mu\in\C$. This gives us nontrivial intertwining operators
of types \fusion{\Mp}{\M^{\pm}}{\M^{\pm}},
\fusion{\M^{\pm}}{\Ml}{\Ml} and \fusion{\Ml}{\Mu}{M(1,\lm+\mu)}. The fusion rules
of corresponding types are nonzero. In addition, in
\cite{FLM},  twisted vertex operator from $\Ml$ to $\Hom(\Mt,\Mt)\{z\}$
was obtained for every $\lm\in\C$. This provide the nontrivial
intertwining operators of types \fusion{\M^{\pm}}{\Mt^{\pm}}{\Mp}
and \fusion{\Ml}{\Mt^{\alpha}}{\Mt^{\beta}} for any
$\alpha,\beta\in\{+,-\}$. Thus the fusion rules of corresponding types
are also nonzero. The study of the contractions also shows that all nonzero
fusion rules are derived from these fusion rules by means of the above
symmetry of fusion rules and the equivalency between $\Ml$ and
$M(1,-\lm)$ for $\lm\in\C$.

The organization of this paper is as follows: We recall definitions of vertex operator
algebras, modules and fusion rules in Subsection \ref{S2.1}, those of Zhu's algebras and
its bimodule in Subsection \ref{S2.2}, where we also explain the relation between
fusion rules and the bimodules, and  we review the vertex
operator algebra
$\Mp$ and its irreducible modules in Subsection \ref{S3}. In Subsection \ref{S4.1}, we
describe the irreducible decompositions of the irreducible
$\Mp$-modules as modules for Virasoro algebra and prove that
$N_{M^{1}M^{2}}^{M^{3}} $ is invariant under the any permutations
of $\{1,2,3\}$ for irreducible $\Mp$-modules $M^{i}$ ($1\leq i\leq
3$), in Subsection \ref{S2.3},  we prove some lemmas and proposition
(Proposition \ref{P3.5}) which gives a generalization of Zhu's anti-isomorphism of
$A(V)$, and in Subsection
\ref{S4.2}, we give a set of generators of
$A(M)$ for irreducible $\Mp$-module $M$, and show that all fusion rules among
irreducible $\Mp$-modules are less than 2. In Subsection \ref{S5.1}, we explain
that the vertex operators and twisted vertex operators constructed in
\cite{FLM}  give some nonzero intertwining operators among irreducible
$\Mp$-modules, and state the main theorem. In Subsection \ref{S5.2}, we prove
the main theorem by studying the structure of the contraction
$(M_{0}^{3})^{*}\cdotp A(M^{1})\cdotp M_{0}^{2}$ for irreducible
$\Mp$-modules $M^{i}$ ($1\leq i\leq 3$).


\section{Preliminaries}\label{S2}

We recall the definitions of vertex operator algebras, its modules
from \cite{FLM}, \cite{DLM1} and \cite{DMZ},  and fusion rules  from
\cite{FHL} in Subsection
\ref{S2.1}. In Subsection \ref{S2.2}, following \cite{Z} and \cite{FZ}, we review
the definition of Zhu's algebra
$A(V)$ associated to a vertex operator algebra $V$ and its bimodule $A(M)$  
for an $\N$-gradable $V$-module $M$. In Subsection \ref{S3}, following
\cite{FLM}, we recall the vertex operator algebra $\Mp$ and its irreducible modules.

Throughout the paper, $\N$ is the set of
nonnegative integers and $\Z_{>0}$ is the set of positive integers.
For vector space $V$, the vector space of formal power
series in $z$ is denoted by
\eqa 
V\{z\}=\left\{\left. \,\sum_{n\in\sC}
v_{n}z^{n}\,\right|\, v_{n}\in V\,\right\},\nn\eeqa
and we set the subspaces $V[[z,z^{-1}]]$ and $V((z))$ as follows:
$$ 
V[[z,z^{-1}]]=\left\{\left.\,\sum_{n\in\sZ}
v_{n}z^{n}\,\right|\,v_{n}\in V\,\right\},\ V((z))=\left\{\left.\,
\sum_{n=k}^{\infty} v_{n}z^{n}\,\right|\,k\in\Z, v_{n}\in V\,\right\}.$$
For $f(z)=\sum_{n\in\sC}v_{n}z^{n}\in V\{z\}$, $v_{-1}$ is called the formal residue
denoted by $\Res_{z} f(z)$$=v_{-1}$. 

\subsection{Vertex operator algebras, modules and fusion rules}\label{S2.1}
\bd{D2.0}A $\Z$-graded vector space $V=\Ops_{n\in \sZ} V_{n}$ such that
$\dim V_{n}$ is finite for all integer $n$ and $V_{n}=0$ for
sufficiently small integer $n$ is called a {\rm vertex operator algebra}
if  $V$ is equipped with a linear map 
\eqa
Y:V& \to & (\End V)[[z, z^{-1}]] \nn \\
v & \mapsto & Y(v,z)=\sum_{n\in \sZ} v_{n} z^{-n-1}\ (v_{n}\in \End
V)\nonumber
\eeqa
and with two distinguished vectors $\1\in V_{0}$ and $\ome\in V_{2}$ such that the
following conditions hold for
$a,\ b\in V\ and\ m,\ n\in\Z$:
\eqa
Y(a,z)b\in V((z)),\nn
\eeqa

\noindent
$(${\rm Jacobi identity}$)$
\eqa
&&z_{0}^{-1}\delta\left( {z_{1}-z_{2}\over
z_{0}}\right)Y(a,z_{1})Y(b,z_{2})-z_{0}^{-1}\delta\left(
{z_{2}-z_{1}\over -z_{0}}\right)Y(b,z_{2})Y(a,z_{1})\nn\\
&&{}=z_{2}^{-1}\delta\left( {z_{1}-z_{0}\over
z_{2}}\right)Y(Y(a,z_{0})b,z_{2}),\nn\eeqa
where $\delta(z)=\sum_{n\in\sZ} z^{n}$ and all binomial expressions are to be
expanded as formal power series in the second variable.
\eqa
Y(\1,z)=\id _{V},\ Y(a,z)\1\in V[[z]] \ {\it and} \ Y(a,z)\1|_{z=0}=a.\nn\eeq

We set $Y(\ome, z)=\sum_{n\in \sZ} L(n)\ z^{-n-2}$, then
$L(n)$, $(n\in \Z) $ form a Virasoro algebra,
\eqa
[L(m), L(n)]=(m-n)L(m+n)+{m^{3}-m\over 12}\delta_{m+n,0}
c_{V}\label{2.6}
\eeqa
for any $m,\ n\in \Z$, where $c_{V}\in\C$ which is called the
{\it central charge} of $V$.
\eqa
L(0)a=na\ \ {\it for}\ n\in \Z,\ a\in V_{n},\nn\\
Y(L(-1)a,z)={d\over dz} Y(a,z).\label{2.8}
\eeqa
\ed

The vertex operator algebra is denoted by $(V,\ Y,\ \1,\ \ome)$ or
simply by $V$. An element $a\in V_{n}$ is called a {\it
homogeneous element of weight} $n$ denoted by $n=\wt(a)$.
 
An {\em automorphism}
$g$ of a vertex operator algebra $V$ is a linear automorphism of $V$ such
that  $g\, Y(a, z)\, g^{-1}=Y(g(a),z)$ for all $a\in V$ and $g(\ome)= \ome$.
Set $\Aut V$ to be the set of all automorphisms of $V$ and let $G$ be a
subgroup of $\Aut V$. Then the fixed point set for $G$
naturally becomes a vertex operator algebra. This vertex operator
algebra is called the {\it orbifold} of $V$ (cf. \cite{DVVV}, \cite{DM}).

Let $g$ be an
automorphism of a vertex operator algebra $V$ of order $T$. Then $V$
is decomposed into the eigenspaces for $g$:
$$
V=\Ops_{r=0}^{T-1} V^{r},\ V^{r}=\{\, a\in V\, |\, g(a)=e^{-{2 \pi i\over
T}} a\, \}.$$
\bd{D2.2}\ Let $V$ be a vertex operator algebra and let $g$ be an
automorphism of order $T$. A {\em weak $g$-twisted
$V$-module} $M$ is a vector space equipped with a linear map
\eqa
Y_{M}:V&\to&(\End M)\{z\},\nn\\
a &\mapsto&Y_{M}(a,z)=\sum_{n\in\sQ} a_{n}^{M}z^{-n-1},\
(a_{n}^{M}\in \End M)\nn
\eeqa
such that the following conditions hold for $0\leq r\leq T-1,\ a\in
V^{r},\ b\in V\ and\ u\in M$:
\eqa
Y_{M}(a,z)=\sum_{n\in {r/ T}+{\sZ}} a_{n}^{M} z^{-n-1},\
Y_{M}(a,z)v\in z^{-{r\over T}}M((z)),\ {\it and}\ Y_{M}(\1,z)=\id _{M}.\nn
\eeqa
\noindent
$($\rm twisted Jacobi identity$)$
\eqa
&&z_{0}^{-1}\delta\left( {z_{1}-z_{2}\over
z_{0}}\right)Y_{M}(a,z_{1})Y_{M}(b,z_{2})-z_{0}^{-1}\delta\left(
{z_{2}-z_{1}\over -z_{0}}\right)Y_{M}(b,z_{2})Y_{M}(a,z_{1})\nn\\
&&{}=z_{2}^{-1}\left( {z_{1}-z_{0}\over
z_{2}}\right)^{-{r\over T}}\delta\left( {z_{1}-z_{0}\over
z_{2}}\right)Y_{M}(Y(a,z_{0})b,z_{2}),\label{2.17}
\eeqa
\ed
The weak $g$-twisted $V$-module is denoted by $(M,Y_{M})$ or
simply by $M$. In the case $g$ is identity of $V$, the weak $g$-twisted
$V$-module is called a {\it weak
$V$-module}. Here and further we write the component operator
$a_{n}^{M}\ (a\in V,\ n\in\Q)$ by $a_{n}$ for simplicity.
 
Let $(M,Y_{M})$ be a weak $g$-twisted $V$-module and set
$Y_{M}(\ome ,z)=\sum_{n\in\sZ} L(n)z^{-n-2}$. Then the operators
$\{\,L(n); (n\in \Z),\id_{M}\,\}$ also form Virasoro algebra with central charge $c_{V}$.
Moreover we also have (\ref{2.8}). (see \cite{DLM2}).

\bd{D2.4} Let $V$ be a vertex operator algebra and let $g$ be an
automorphism of $V$ of order $T$. A {\em ${1\over
T}\N$-gradable $g$-twisted $V$-module} $M$ is a weak 
$g$-twisted $V$-module which has a ${1\over
T}\N$-grading $M=\ops_{n\in {1\over
T}\ss{\sN}} M_{n}$ such that $
a_{m}M_{n}\subset M_{{\ss\wt}(a)+n-m-1}
$
holds  for any homogeneous $a\in V,\ n\in{1\over T}{\N}\ and\ m\in\Q$.
\ed

In the case $g$ is identity on $V$, the ${1\over
T}\N$-gradable $g$-twisted $V$-module is called an 
{\em $\N$-gradable $V$-module}. An element $u\in
M_{n}$ is called a {\it homogeneous element of degree} $n$ denoted by $n=\deg(v)$.

\bd{D2.5}  An {\em  ordinary $g$-twisted
$V$-module} $M$  is a weak $g$-twisted $V$-module on which $L(0)$ acts
semisimply:
\eqa
M=\Ops_{\lm \in \sC} M(\lm),\ M(\lm)=\{\, u\in M\,|\, L(0) u=\lm
u\,\} \nn
\eeqa
such that each eigenspace is finite dimensional and for fixed
$\lm\in\C$, $M(\lm+n/T)=0$ for sufficiently small
integer $n$.
\ed
In the case $g$ is identity of $V$,  the ordinary $g$-twisted
$V$-module is called an {\it ordinary $V$-module} or more simply
a {\it $V$-module}. An element $u\in M(\lm)$ is said to be a {\it homogeneous
of weight} $\lm$ denoted by $\lm = \wt(u)$.

The notions of submodules, irreducible modules are defined in the obvious way. Let
$M$ be a $V$-module, then the restricted dual $M'=\ops_{\lm\in\C} M(\lm)^{*}$ is
 a $V$-module, and the vertex operator $Y_{M}^{*}(a,z)$ for $a\in V$ is defined by
$$\B Y_{M}^{*}(a,z)u',v\K =\B u',Y_{M}(e^{z L(1)} (-z^{-2})^{L(0)}a,z^{-1})v\K$$ for
$u'\in M',\ v\in M$. This $V$-module $(M',Y_{M}^{*})$ is called the {\em contragredient
module} of $M$. It is known that if $M$ irreducible, then $M'$ is also irreducible (cf.
\cite{FHL}). 

\bd{D3.7} Let $V$ be a vertex operator algebra and $(M^{i},
Y_{M^{i}})\ (i=1,2,3)$ be weak $V$-modules. An
{\em intertwining operator of type} \fusion{M^{1}}{M^{2}}{M^{3}}
 is a linear map
\eqa
I:M^{1}&\to&(\Hom (M^{2},M^{3}))\{z\},\nn\\
v\ \ &\mapsto&I(v,z)=\sum_{n\in\sC}v_{n} z^{-n-1}\
(v_{n}\in\Hom(M^{2},M^{3}))\nn
\eeqa
such that for $a\in V, v\in M^{1}\ and\ u\in M^{2}$, following
conditions hold:

\noindent
For fixed $n\in\C$, $v_{n+k}u=0$ for sufficiently large integer $k$,

$($Jacobi\ identity$)$
\eqa
&&z_{0}^{-1}\delta\left( {z_{1}-z_{2}\over
z_{0}}\right)Y_{M^{3}}(a,z_{1})I(v,z_{2})-z_{0}^{-1}\delta\left(
{z_{2}-z_{1}\over -z_{0}}\right)I(v,z_{2})Y_{M^{2}}(a,z_{1})\nn\\
&&{}=z_{2}^{-1}\delta\left( {z_{1}-z_{0}\over
z_{2}}\right)I(Y_{M^{1}}(a,z_{0})v,z_{2}),\label{3.16}
\eeqa
$$\displaystyle {d\over dz} I(v,z)= I(L(-1)v,z).$$
\ed
We denote the vector space consisted of the intertwining operators of
type \fusion{M^{1}}{M^{2}}{M^{3}} by
$I$\fusion{M^{1}}{M^{2}}{M^{3}}. The dimension of this vector space
 is called a {\it
fusion rule} of corresponding type denoted by
$N_{M^{1}M^{2}}^{M^{3}}$. It is well known that fusion rules have the following
symmetry (see \cite{FHL} and \cite{HL}).
\bp{P3.13} Let $M^{i}$ $(i=1,2,3)$ be $V$-modules. Then 
$$N_{M^{1}M^{2}}^{M^{3}}=N_{M^{2}M^{1}} ^{M^{3}}\ \ and\ \
N_{M^{1}M^{2}}^{M^{3}} = N_{M^{1}(M^{3})'}^{(M^{2})'}.$$
\ep

\subsection{Zhu's algebra $A(V)$ and $A(V)$-bimodule $A(M)$}\label{S2.2}
We recall the definition of Zhu's algebra. Two bilinear products $*$ and $\circ$
on
$V$ are defined as follows: For homogeneous $a\in V$, and $b\in V$, we define 
$$
a*b=\left( \Res_{z} {(1+z)^{{\swt}(a)}\over z}
Y(a,z)\right)b,\ \ 
a\circ b=\left( \Res_{z} {(1+z)^{{\swt}(a)}\over z^{2}}
Y(a,z)\right)b$$
and extend to $V$ by linearity. Let $O(V)$ be the linear
span of $a\circ b$ ($a,b\in V$) and set
$A(V)=V/O(V)$.
Let $M$ be an $\N$-gradable $V$-module. For every homogeneous $a\in V$,
define $ o(a) =a^{M}_{{\swt}(a)-1}$ and extend to $V$ linearly. 
The following proposition is due to Zhu (see \cite{Z}).
\bp{P3.1} (1) The bilinear product $*$ induces $A(V)$ an associative
algebra structure. The vector $\1+O(V)$ is the identity and $\ome
+O(V)$ is in the center of $A(V)$.

\noindent
(2) The linear map $o:V\rightarrow \End M_{0}$, $v\mapsto o(v)|_{M_{0}}$ induces
an associative algebras homomorphism $o:A(V)\to \End M_{0}$. Thus $M_{0}$ is 
a left $A(V)$-module.
\ep

Next we recall $A(V)$-bimodule $A(M)$. Let $(M,Y_{M})$ be an $\N$-gradable
$V$-module. Define bilinear maps $
*:V\ts M\to M,\ \ \circ:V\ts M\to M,\ {*}:M\ts V\to M\nn$ by
\eqa 
a*u&=&\left( \Res_{z} {(1+z)^{{\swt}(a)}\over z}
Y_{M}(a,z)\right)u=\sum_{i=0}^{{\swt}(a)} {\wt (a)\choose i}
a_{i-1}u,\label{3.4}\\
a\circ u&=&\left( \Res_{z} {(1+z)^{{\swt}(a)}\over z^{2}}
Y_{M}(a,z)\right)u=\sum_{i=0}^{{\swt}(a)} {\wt (a)\choose i}
a_{i-2}u,\nn\\
u*a&=&\left( \Res_{z} {(1+z)^{{\swt}(a)-1}\over z}
Y_{M}(a,z)\right)u=\sum_{i=0}^{{\swt}(a)-1} {\wt (a)-1\choose i}
a_{i-1}u,\label{3.6}\eeqa
for homogeneous $a\in V$ and $u\in M$ respectively, and extend
to linear operations on $V$. Let $O(M)$ be the linear span of $a\circ u$ for $a\in
V$, $u\in M$, and set $A(M)=M/O(M)$. We denote the image of $u\in M$ in $A(M)$ by
$[u]$. Then $A(M)$ is
$A(V)$-bimodule, and the left and right actions are given by $[a]*[u]=[a*u]$ and
$[u*a]=[u*a]$ respectively for $a\in V, u\in M$.
\bd{D3.10} Let $A$ be an associative algebra and $R,B,L$ be a right
$A$-module, an $A$-bimodule and a left $A$-module respectively.
The tensor product of $R,B$ and $L$ as $A$-module $R\ots_{A} B\ots_{A}L$ is called a
{\em contraction} of $R,B$ and $L$ and denote by $ R\cdotp B\cdotp L$
\ed

Now let $M^{i}=\ops_{n=0}^{\infty} M_{n}^{i}\ (i=1,2,3)$ be
$\N$-gradable $V$-modules. Suppose that for any $i$, there exists a scalar $h_{i}\in\C$
such that
$L(0)$ acts on
$M_{n}^{i}$ as $h_{i}+n$. Let $I$ be an intertwining
operator of type \fusion{M^{1}}{M^{2}}{M^{3}}. Then for each $v\in M^{1}$, we have 
$I(v,z)\in z^{-h_{1}-h_{2}+h_{3}}(\Hom (M^{2},M^{3}))[[z,z^{-1}]]$ (see Proposition 1.5.1
 in
\cite{FZ}). We define
$o(v)=\Res_{z} z^{h_{1}+h_{2}-h_{3}+\deg(v)-1}I(v,z)$ for homogeneous $v\in
M^{1}$ and extend to $M^{1}$ by linearity. Then we have following theorem (see
Theorem 1.5.2 in \cite{FZ}).
\bt{T3.11}  Let $M^{i}=\ops_{n=0}^{\infty} M_{n}^{i}\ (i=1,2,3)$ be
$\N$-gradable $V$-modules. Suppose that for each $M^{i}$, 
there exists an $h_{i}\in\C$ such that $L(0)$ acts on
$M_{n}^{i}$ as a scalar $h_{i}+n$. Let  
$(M^{3})'=\ops_{n=0}^{\infty} (M_{n}^{3})^{*}$ be the contragredient
module of $M^{3}$. Then a linear map 
\eqa 
\pi:I{\fusion{M^{1}}{M^{2}}{M^{3}}}\to ((M^{3}_{0})^{*}\cdotp A(M^{1})\cdotp
M^{2}_{0})^{*},\ I\mapsto\pi(I).\label{3.21}
\eeqa
is defined by $\pi(I)(v_{3}'\ots [v_{1}]\ots v_{2})=\B
v_{3}',o(v_{1})v_{2}\K $ for $I\in I\fusion{M^{1}}{M^{2}}{M^{3}},\ v_{3}'\in
(M^{3}_{0})^{*},\ v_{1}\in M^{1}$ and
$v_{2}\in M^{2}_{0}$.
\et
  
If $V$ is rational, that is, all $\N$-gradable $V$-modules are completely reducible,
then the linear map $\pi$ is an isomorphism. In general, $\pi$ is not surjective (see
\cite{L}), but we have the following proposition. 
\bp{P3.12} Let $M^{i}=\ops_{n=0}^{\infty} M_{n}^{i}\ (i=1,2,3)$ be as in Theorem
$\ref{T3.11}$. Suppose that $M^{2}$  is an irreducible and that $M^{3}$ is an irreducible
ordinary. Then
$\pi$ is injective. Thus we have 
\eqa N_{M^{1}M^{2}}^{M^{3}}\leq \dim ((M^{3}_{0})^{*}\cdotp 
A(M^{1})\cdotp M^{2}_{0})^{*}.\nn
\eeqa
\ep
\proof Let $I$ be an intertwining operator of type \fusion{M^{1}}{M^{2}}{M^{3}}.
Suppose that $\pi (I)=0$. Then we have 
$\B v_{3}',I(v_{1},z)v_{2}\K =0$  for any $v_{3}'\in  (M^{3}_{0})^{*},\ v_{1}\in M^{1}$
and $v_{2}\in M^{2}_{0}$.  Set
 $$W=\{u\in M^{2}| \B v_{3}',I(v_{1},z)u\K=0\  {\rm for\ any}\ v_{3}'\in
(M^{3}_{0})^{*}, v_{1}\in M^{1}\}\supset M^{2}_{0}.$$
and fix a nonzero $u\in M_{0}^{2}$. Then we have $a_{n}u\in W$ for any homogeneous
$a\in V$, $n\in \Z$. In fact, it is obvious that $a_{n}u\in W$ for $n\geq \wt(a)-1$,
and   for $n<\wt(a)-1$,\ we have $\B v_{3}, a_{n}(I(v_{1},z)u)\K=0$ for
any $v_{3}'\in (M^{3}_{0})^{*}$ and $v_{1}\in M^{1}$. Hence we see
that 
$$\B v_{3}',I(v_{1},z)a_{n}u\K=\B
v_{3}',a_{n}I(v_{1},z)u\K-\sum_{i=0}^{\infty}{n\choose i}
z^{n-i}\B v_{3}', I(a_{i}v_{1})u\K=0,$$
from Jacobi identity (\ref{3.16}). Thus $a_{n}u\in W$ for all $n<\wt(a)-1$. 
Since $M^{2}$ is irreducible, $M^{2}$ is spanned by $a_{n}u$ for homogeneous $a\in
V$ (see \cite{DM}), and then $W=M^{2}$. 

Let $w'\in (M^{3}_{0})^{*}$ be a nonzero element. Then for every $a\in V$ and $n\in
\Z$, we have 
$$ \B w', a_{n}(I(v_{1},z)v_{2})\K=\B w', I(v_{1},z)a_{n}v_{2}\K+
\sum_{i=0}^{\infty}{n\choose i}z^{n-i}\B w',
I(a_{i}v_{1},z)v_{2}\K=0.$$
Hence $\B w', Y_{M^{3}}(a,z_{0})I(v_{1},z)v_{2}\K=0$. This implies that 
\eqa 
\B Y_{M^{3}}^{*}(a,z_{0})w',I(v_{1},z)v_{2}\K=\B w',
Y_{M^{3}}(e^{z_{0}L(-1)}(-z_{0}^{-2})^{L(0)}a,z_{0}^{-1})I(v_{1},z)v_{2}
\K=0.\nn
\eeqa
Therefore $\B a^{*}_{n}w',I(v_{1},z)v_{2}\K=0$ for any $a\in V$ and
$n\in\Z$, where $Y_{M^{3}}^{*}(a,z)=\sum_{n\in\sZ}a_{N}^{*}z^{-n-1}$. Since $(M^{3})'$ 
is irreducible, it is spanned by $a^{*}_{n}w'$ for $a\in V$ and
$n\in\Z$. Then for every $v'_{3}\in( M^{3})'$, $v_{1}\in M^{1}$ and
$v_{2}\in M^{2}$,  $\B v'_{3},I(v_{1},z)v_{2}\K=0$, which means $I=0$.\qed

\subsection{Vertex operator algebra $\Mp$}\label{S3}

Let $\h$ be a $d$-dimensional vector space with a
nondegenerate symmetric bilinear form $\B\,,\K$, and
let $\hh=\h\ots\C[t,t^{-1}]\ops \C K$ be a Lie algebra with the commutation relation
given by $ [h\ots t^{m}, h'\ots t^{n}]=m\delta_{m+n,0}\B h,h'
\K K$, $[K,\hh]=0$. Set $\hp=\h\ots\C[t]\ops\C K$. For $\lm\in\h$, let $\C e^{\lm}$
be a one-dimensional $\hp$-module on which $\h\ots t\C[t]$ acts trivially, $\h$ acts
as $\B h,\lm\K$ for $h\in\h$ and $K$ acts as $1$. Let $M(1,\lm)$ be an $\hh$-module
induced by the
$\hp$-module 
$\C e^{\lm}$:
\eqa M(1,\lm)=U(\hh)\ots_{U({\ss{\shp}})}\C e^{\lm}\simeq S(\h\ots
t^{-1}\C[t^{-1}])\ ({\rm linearly}).\nn
\eeqa

Denote the action of $ h\ots t^{n}$, $(h\in\h,n\in \Z)$ on $\Ml$ by $h(n)$ and
set  $h(z)=\sum_{n\in\sZ} h(n) z^{-n-1}$.
For $\lm,\mu\in\h$, we define a linear map
$P_{\lm\mu}:\Mu\to M(1,\lm+\mu)$ by
$P_{\lm\mu}(u\ots e^{\mu})=u\ots e^{\lm+\mu}$ for $u\in S(\h\ots t^{-1}\C[t^{-1}])$.
Then a {\it vertex operator} associated with $e^{\lm}$ is defined by
\eqa
I_{\lm\mu}(e^{\lm},z)=\exp\left(\sum_{n\in
\sZ_{>0}}{\lm(-n)\over n} z^{n}\right)\exp\left(-\sum_{n\in
\sZ_{>0}}{\lm(n)\over n} z^{-n}\right)P_{\lm\mu}z^{\B\lm,\mu\K}.\label{4.3}
\eeqa
The vertex operator associated with $h_{1}(-n_{1})h_{2}(-n_{2})\cdots h_{k}(-n_{k})
e^{\lm}\in\Ml$
$(n_{i}\in\Z_{>0},\ h_{i}\in\h)$ is defined by
\eqa
&&I_{\lm\mu}(h_{1}(-n_{1})h_{2}(-n_{2})\cdots h_{k}(-n_{k})
e^{\lm},z)\nn\\
&&{}=\NO\der^{(n_{1}-1)}h_{1}(z)\der^{(n_{2}-1)}h_{2}(z)
\cdots\der^{(n_{k}-1)}h_{k}(z)Y_{\lm}(e^{\lm},z)\NO,\label{4.4}
\eeqa
where $\der^{(n)}={1\over n!}({d\over dz})^{n}$ for $n\in
\N$ and the {\it normal ordering} $\NO\cdotp\NO$ is an
operation which reorders so that $h(n)$ $(n<0)$, $P_{\lm}$ to
be placed to the left of $h(n)$ $(n\in\N)$. We extend $I_{\lm\mu}$ to $\Ml$ by
linearity. Now let $\{\al_{1},\al_{2},\ldots,\al_{d}\}$ ($d=\dim \h$) be an
orthonormal basis of $\h$, and set
$\1=1\ots e^{0}$, $\ome=(1/2)\sum_{i=1}^{d}\al_{i}(-1)^{2}\1
\in M(1,0)$.
Then $(M(1,0),I_{00},\1,\ome)$ is a simple
vertex operator algebra with central charge $d$, and
$\{(\Ml,I_{0\lm})\,|\,\lm\in\h \}$ gives all irreducible 
$M(1,0)$-modules (see \cite{FLM}). The vertex operator algebra $M(1,0)$ is called
the {\it free bosonic vertex operator algebra}, and denoted by \M.

Let $\th$ be an automorphism of $\M$ defined by 
$$\th(h_{1}(-n_{1})h_{2}(-n_{2})\cdots h_{k}(-n_{k})\1)=
(-1)^{k}h_{1}(-n_{1})h_{2}(-n_{2})\cdots
h_{k}(-n_{k})\1$$
for $h_{i}\in \h,\ n_{i}\in \Z_{>0}$. We denote the orbifold of
$\M$ for $\th$ by $\Mp$ and the $-1$-eigenspace of $\M$ by $\Mm$. Then we
have following  proposition.
\bp{P4.3} The vertex operator algebra
$\Mp$ is simple, and $\Mm$, $\Ml$ $(\lm\neq0)$ are irreducible
$\Mp$-modules. Moreover
$M(1,\lm)$ is isomorphic to $M(1,-\lm)$ as $\Mp$-module.
\ep
\proof This is proved by using Theorem 2, Theorem 4.4 and
Theorem 6.1 of \cite{DM}.\qed
\vs

Next we consider the $\th$-twisted $\M$-module. Let $\hm= \h\ots t^{{1\over
2}}\C[t,t^{-1}]\ops\C K$ be a Lie algebra with commutation relation $[h\ots t^{m},
h'\ots t^{n}]=m\delta_{m+n,0}\B h,h' \K K$, $[K,\hm]=0$ for $h,h'\in \h,\ m,n\in
1/2+\Z$. Set
$\hmp=\h\ots t^{1/2}\C[t]\ops\C K$. Then $\C$ is viewed as an
$\hmp$-module on which $\h\ots t^{1/2}\C[t]$ acts trivially and $K$ acts as $1$.
Set $\Mt$ to be the induced $\hm$-module:
$$\Mt=U(\hm)\ots_{U({{\shmp}})}\C\simeq S(\h\ots
t^{-{1\over 2}}\C[t^{-{1\over 2}}])\ ({\rm linearly}).$$

Denote the action of $h\ots t^{n}$ $(h\in\h,\ n\in {1/2}+\Z)$ on
$\Mt$ by $h(n)$ and set $h(z)=\sum_{n\in {1\over 2}+\sZ}
h(n)z^{-n-1}$. For $\lm\in\h$, a {\it twisted vertex
operator} associated with $e^{\lm}\in\Ml$ is defined as follows:
\eqa I_{\lm}^{\th}(e^{\lm},z)=z^{-{\B\lm,\lm\K\over 2}} 
\exp\left(\sum_{n\in{1\over 2}+\sN}{\lm(-n)\over
n}z^{n}\right)
\exp\left(-\sum_{n\in{1\over 2}+\sN}{\lm(n)\over n}
z^{-n}\right)\label{4.125}
\eeqa
For $h_{1}(-n_{1})h_{2}(-n_{2})\cdots
h_{k}(-n_{k})\ots e^{\lm}\in\Ml$ $(h_{i}\in\h$, $n\in \Z_{>0})$, set
\eqa
&&W_{\lm}^{\th}(h_{1}(-n_{1})h_{2}(-n_{2})\cdots h_{k}(-n_{k})
e^{\lm},z)\nn\\
&&{}=\NO\der^{(n_{1}-1)}h_{1}(z)\der^{(n_{2}-1)}h_{2}(z)
\cdots\der^{(n_{k}-1)}h_{k}(z)I_{\lm}(e^{\lm},z)\NO,\nn
\eeqa
and extend to a linear operator on $\Ml$, where the normal
ordering $\NO\cdotp\NO$ is an operation which sifts $h(n)$
$(n\in\Z_{>0})$ to right and $h(n)$ $(n\in\Z_{<0})$ to left.
Let $c_{mn}\in\Q$ $(m,n\in\N)$ be constants defined by
the following formal power series expansion
\eqa \sum_{m,n\geq0} c_{mn}x^{m}y^{n}=-\log \left(
{(1+x)^{1\over 2}+(1+y)^{1\over 2}\over 2}\right),\nn
\eeqa
and define an operator $\Del_{z}$ on $\Ml$ by $
\Del_{z}=\sum_{m,n\geq 0}\sum_{i=0}^{d} c_{mn}
\al_{i}(m)\al_{i}(n) z^{-m-n}.$
Then the twisted vertex
operator associated with $u\in \Ml$ is defined by
\eqa
I_{\lm}^{\th}(u,z)=W_{\lm}^{\th}(e^{\Del_{z}}u,z).\label{4.13}\eeqa

Now define the action of $\th$ on $\Mt$ by
\eqa \th(h_{1}(-n_{1})h_{2}(-n_{2})\cdots h_{k}(-n_{k})\1)=
(-1)^{k}h_{1}(-n_{1})h_{2}(-n_{2})\cdots
h_{k}(-n_{k})\1\label{4.16}
\eeqa
for $h_{i}\in \h,\ n_{i}\in 1/2+\N$, and denote the $\pm
1$-eigenspace of $\Mt$ by $\Mt^{\pm}$ respectively.
Then we have following proposition.
\bp{P4.5}(1) $(\Mt, I_{0}^{\th})$ is an irreducible $\th$-twisted $\M$-module.\\ 
(2) $(\Mt^{\pm},I_{0}^{\th})$ are irreducible $\Mp$-modules.
\ep
\proof  Part (1) is one of results of Chapter 9 in \cite{FLM}, and  part (2) is a direct
consequence of Theorem 5.5 of \cite{DLi}.\qed
\vs

From now on, we consider the case $\dim \h=1$ and set
$\al_{1}=h$. For $\lm\in\C$, we denote the irreducible
$\M$-module $M(1,\lm h)$ by $\Ml$ and $e^{\lm h}$ by $e^{\lm}$. 

Set $J=h(-1)^{4}\1-2 h(-3)h(-1)\1+(3/2) h(-2)^{2}\1\in \Mp$ which
is lowest weight vector of weight $4$ for Virasoro algebra. Then the vertex operator
algebra $\Mp$ is generated by
$\ome $ and $J$ (see \cite{DG}), and the Zhu's algebra $A(\Mp)$ is generated by
$[\ome]$ and $[J]$. More precisely, there is an isomorphism  $
\C[x,y]/I\stackrel{\sim}{\to} A(\Mp),\ x+I\mapsto[\ome],\ y+I\mapsto[J]$, where $I$
is an ideal generated by two polynomials $(y-4 x^{2}+ x)(70 y+908 x^{2}-515 x+27)$,
$(y-4 x^{2}+x)(x-1)(x-{1/16})(x-{9/16})$ (see \cite[Theorem 4.4]{DN}).

In \cite[Theorem 4.5]{DN}, all equivalence classes of
irreducible $\N$-gradable $\Mp$-modules are classified as follows: 
\bt{T4.8} The set $\{ \M^{\pm}, \Mt^{\pm},
\Ml(\simeq M(1,-\lm));(\lm\neq0)\}$ gives all inequivalent irreducible
$\Mp$-modules.
\et

Let $M$ be an irreducible $\Mp$-module. Then the homogeneous space of
the lowest weight (written by $a_{M}\in\C$)  is one dimensional, and
$M$ has naturally an $\N$-gradation defined by $M_{n}=M(a_{M}+n)$
for any $n\in\N$, where $M(\lm)$ is the homogeneous space
of weight $\lm$ of $M$. One of the basis (written by
$v_{M}$) of $M_{0}$ is given by Table \ref{Table4.1}, and we denote the
dual basis of $(M_{0})^{*}$ corresponding to $v_{M}$ by $v_{M}'$. Since
$M_{0}$ is one-dimensional, 
$[J]$ also acts on this space as scalar denoted by $b_{M}\in\C$. 
\begin{table}[htbp]
\begin{center}
\begin{tabular}{|c|c|c|c|c|c|}\hline
$M$&$\Mp$&$\Mm$&$\Ml \ (\lm\neq 0)$&$\Mtp$&$\Mtm$\\
\hline
 $v_{M}$&$\1$&$h(-1)\1$&$e^{\lm}$&$1$&$h(-1/2)1$\\ \hline
$a_{M}$&$0$&$1$&$\lm ^{2}/2$&$1/16$&$9/16$\\ \hline
$b_{M}$&$0$&$-6$&$\lm^{4}-\lm^{2}/2$&$3/128$&$-45/128$\\
\hline
\end{tabular}
\end{center}
\caption[]{The basis of top levels and the eigenvalues of $[\ome]$ and
$[J]$.}
\label{Table4.1}
\end{table}


\section{A spanning set of $A(M)$}\label{S4}
This section is divided into three subsections. In first subsection, we describe the
irreducible decompositions of irreducible $\Mp$-modules as modules for Virasoro
algebra. In second subsection, we prove some lemmas. In last subsection, we give a
spanning set of $A(M)$ for irreducible
$\Mp$-module $M$.

\subsection{Irreducible decompositions of irreducible $\Mp$-modules}\label{S4.1}
Let $W$ be a module for Virasoro algebra with central charge $c\in\C$
such that $L(0)$ acts on $W$ semisimply and each eigenspace is
finite dimensional. Then  the {\it formal character} of $W$ is defined by 
\eqa 
\ch_{W} ={\rm tr}_{W}\, q^{L(0)-{c\over 24}}=q^{-{c\over
24}}\sum_{\lm\in\C}(\dim W(\lm)) q^{\lm}\in \Z\{z\},\nn
\eeqa
where $W(\lm)$ is the eigenspace of weight $\lm$ for $L(0)$. Let
$L(1,\lm)$ be the irreducible lowest weight module for Virasoro
algebra with central charge $1$ and lowest weight $\lm$. Then the formal 
character of $L(1,\lm)$ is given by
\eqa\ch_{L(1,\lm)}=\left\{
\begin{array}{ll}
{q^{\lm}\over \eta (q)}& \mbox{${\rm if}\ \lm\neq
{n^{2}\over 4}\ {\rm for\ any}\ n\in\Z$} \\
{1\over \eta(q)}\left( q^{n^{2}\over 4}-q^{(n+2)^{2}\over
4}\right)&\mbox{${\rm if}\ \lm={n^{2}\over 4}\ {\rm for\ some\
} n\in\N,$}\end{array}\right.\label{4.19}
\eeqa
where $\eta(q)=q^{1/24}\prod _{n=1}^{\infty} (1-q^{n})$ is the
Dedekind $\eta$-function (see \cite{KR}).

It is well known that the irreducible $\M$-module $\Ml$ $(\lm\in\C)$
is completely reducible module for Virasoro algebra, and decomposed into a direct
sum of irreducible modules for Virasoro algebra as follows (see \cite{WY}, \cite{KR}):
\eqa\Ml=\left\{
\begin{array}{ll}
L(1,{\lm^{2}\over 2})& \mbox{${\rm if}\ {\lm^{2}\over 2}\neq
{n^{2}\over 4}\ {\rm for\ any}\ n\in\Z$} \\
\Ops_{p=0}^{\infty} L(1,{(n+2p)^{2}\over 4})&\mbox{${\rm if}\
{\lm^{2}\over 2}={n^{2}\over 4}\ {\rm for\ some\ }
n\in\N.$}\end{array}\right.\label{4.20}
\eeqa
Furthermore in \cite{DG}, it is shown that the irreducible $\Mp$-modules
$\M^{\pm}$ are also completely reducible as modules for Virasoro algebra, and the
irreducible decompositions are given by
\eqa\Mp=\mathop{\Ops}_{p=0}^{\infty} L(1, 4 p^{2}),\ \ 
\Mm=\mathop{\Ops}_{p=0}^{\infty} L(1, (2p+1)^{2}).\label{4.22}
\eeqa

We show that irreducible $\Mp$-modules $\Mt^{\pm}$ are also
completely reducible as modules for Virasoro algebra. To do this, we
consider the character of  $\Mt$.
\bl{L4.9} We have 
\eqa \ch_{\ss{\Mt}}=\prod_{k=1}^{\infty} {q^{{1/16}-{1/24}}\over (1-q^{k-{1\over
2}})}={1\over \eta(q)}
\sum_{p=0}^{\infty}q^{{(2p+1)^{2}/16}}.\label{4.23}
\eeqa
\el
\proof The first equality is clear by the fact that the weight space of
weight $n$ $(n\in 1/16+(1/2)\Z)$ has an basis 
\eqa
\left\{h(-m_{1})\cdots
h(-m_{k})1\left|\begin{array}{l}\mbox{$k\in
\N,\ m_{i}\in {1\over 2}+\N,\ m_{1}\geq
m_{2}\geq\cdots\geq m_{k}>0,$}\\
\mbox{$m_{1}+\cdots+m_{k}=n$}
\end{array}\right.\right\}.\nn
\eeqa

To prove the second equality, we have to show that
\eqa\prod_{k=1}^{\infty}{(1-q^{k})\over
(1-q^{k-{1/2}})}=\sum_{p=0}^{\infty} q^{p(p+1)/4}.\label{4.24}
\eeqa
From  Jacobi's triple product $
\prod_{n=1}^{\infty}(1-q^{2
n})(1+q^{2n-1}z)(1+q^{2n-1}z^{-1})=\sum_{n\in\sZ}q^{n^{2}}
z^{n}$, we see 
$$\prod_{n=1}^{\infty}(1-q^{4n})(1+q^{2n})= \sum_{n=0}^{\infty}q^{n^{2}+n},$$
by substituting $z=q$. 
Replacing $q$ with $q^{1/4}$, we have
\eqa
\prod_{n=1}^{\infty}(1-q^{n}) (1+q^{n/2})=
\sum_{n=0}^{\infty}q^{(n^{2}+n)/4}.\nn 
\eeqa
On the other hand, direct calculations show
$$\prod_{k=1}^{\infty}{(1-q^{k})\over (1-q^{k-{1/
2}})} =\prod_{k=1}^{\infty}(1-q^{k})(1+q^{k/2}).\nn
$$
Thus the equality (\ref{4.24}) holds.\qed
\vs

We define an hermitian form $(\,|\,):\Mt\times\Mt\to\C$
by 
\eqa 
&&(h(-m_{k})^{p_{k}}\cdots
h(-m_{1})^{p_{1}}1|h(-n_{\ell})^{q_{\ell}}\cdots h(-n_{1})^{q_{1}}1)\nn\\
&&\ \ =\left\{
\begin{array}{ll}
\prod_{i=1}^{k}(m_{i})^{p_{i}} p_{i}!& \mbox{${\rm if}\
{k=\ell}\ {\rm and}\ m_{i}=n_{i}, p_{i}=q_{i}\ {\rm for\ all}\ 1\leq i\leq
k$}
\\
 0&\mbox{${\rm otherwise}$}\end{array}\right.\nn
\eeqa
for $k,\ell\in \N$, $m_{i},n_{j}\in
1/2+\N$, $m_{k}\geq\cdots\geq m_{1}$,\
$n_{\ell}\geq\cdots \geq n_{1}$ and
$p_{i},q_{j}\in\N$ $(1\leq i \leq k,1\leq j\leq
\ell)$. Then we can check that the hermitian form $(\,|\,)$ is a
nondegenerate positive-definite contravariant form on
$\Mt$. Thus $\Mt$ is a unitary representation for Virasoro algebra and is completely
reducible as module for Virasoro algebra. Let 
\eqa
\Mt=\Ops_{h\in\sC} L(1,h)^{\ops m_{h}},\ \
m_{h}\in\N\nn
\eeqa
be the irreducible decomposition of $\Mt$ as Virasoro
algebra module, where 
$$ L(1,h)^{\ops m_{h}}=\overbrace{L(1,h)\ops\cdots\ops
L(1,h)}^{m_{h}}.$$
Since weights which appear in $\Mt$ are in $1/16+(1/2)\Z$, we may
assume that 
$$
\Mt=\Ops_{h\in{1\over 16}+{1\over 2}\sZ} L(1,h)^{\ops
m_{h}},\ \ m_{h}\in\N.$$
In addition, we can show that 
$1/16+n/2\neq m^{2}/4$ for any $n,m\in \Z$. By the formal characters (\ref{4.19}), 
we find that the formal character of $\Mt$ has following form:
$$\ch_{\ss{\Mt}}={1\over \eta(q)}\sum_{h\in{1\over
16}+{1\over 2}\sZ}
 m_{h} q^{h},\ \ m_{h}\in\N.$$
Comparing with  (\ref{4.23}), we have 
$$ m_{h}=\left\{
\begin{array}{ll}
 1 & \mbox{${\rm if}\ h={(2p+1)^{2}\over 16}\ {\rm for\ some}\
p\in\N$}
\\
 0 &\mbox{${\rm otherwise}.$}\end{array}\right.$$
Consequently we have following irreducible decomposition of $\Mt$:
$$
\Mt=\Ops_{p=0}^{\infty} L\left(1,{(2p+1)^{2}\over 16}\right).$$
\bp{P4.10} $\Mp$-modules $\Mt$ and $\Mt^{\pm}$ are completely
reducible as modules for Virasoro algebra and the irreducible
decompositions are given by
\eqa\Mt&=&\Ops_{p=0}^{\infty}
L\left(1,{(2p+1)^{2}\over 16}\right),\label{4.30}\\
\Mtp&=&\Ops_{p=0}^{\infty}\left(
L\left(1,{(8p+1)^{2}\over 16}\right)\Ops L\left(1,{(8p+7)^{2}\over
16}\right)\right),\label{4.31}\\
\Mtm&=&\Ops_{p=0}^{\infty}\left(
L\left(1,{(8p+3)^{2}\over 16}\right)\Ops L\left(1,{(8p+5)^{2}\over
16}\right)\right).\label{4.32}
\eeqa
\ep
\proof We have already proved the irreducible decomposition (\ref{4.30}). By
the definitions of $\Mt^{\pm}$, $\Mtp$ ($\Mtm$) is spanned by all homogeneous
vectors whose weights are in
$1/16+(1/2)\Z$ (resp. $9/16+(1/2)\Z$). Since each weight
which appears in $L(1,(2p+1)^{2}/16)$ is in $(2p+1)^{2}/16+\Z$, we see that if
$(2p+1)^{2}/16$ is in $1/16+(1/2)\Z$ ( $9/16+(1/2)\Z$), then 
$L(1,(2p+1)^{2}/16)$ is contained in $\Mtp$ (resp. $\Mtm$). On one hand, we have 
\eqa {(2p+1)^{2}\over 16}=\left\{
\begin{array}{ll}
 1/16 \pmod{\Z} & \mbox{${\rm if}\
p=0,3\pmod{4}$}
\\
9/16 \pmod{\Z} & \mbox{${\rm
if}\ p=1,2\pmod{4}$}.\nn
\end{array}\right.
\eeqa
Thus we see that 
\eqa \Ops_{p=0}^{\infty}\left(
L\left(1,{(8p+1)^{2}/16}\right)\Ops L\left(1,{(8p+7)^{2}/
16}\right)\right)&\subset&\Mtp,\label{4.33}\\
\Ops_{p=0}^{\infty}\left(
L\left(1,{(8p+3)^{2}/16}\right)\Ops L\left(1,{(8p+5)^{2}/
16}\right)\right)&\subset&\Mtm.\label{4.34}
\eeqa
Since $\Mt=\Mtp\ops\Mtm$, the irreducible decomposition (\ref{4.30}) shows that the
inclusion of (\ref{4.33}) and (\ref{4.34}) are equal. \qed
\vs

From irreducible decompositions (\ref{4.20}), (\ref{4.22}), (\ref{4.31}) and
(\ref{4.32}), we can see that formal characters of irreducible
$\Mp$-modules $\M^{\pm}$, $\Mt^{\pm}$ and $\Ml$ ($\lm\neq0$) are
distinct. So together with Theorem \ref{T4.8}, we
have following lemma.
\bl{L4.11} Let $M$, $N$ be irreducible $\Mp$-modules such that the formal 
character of $M$ coincides with that of $N$. Then $M$\
is isomorphic to $N$. 
\el
In particular, if $M$ is irreducible $\Mp$-module, then the
character of its contragredient module $M'$ is same of $M$. Thus Lemma \ref{L4.11}
shows that $M$ and its contragredient module
$M'$ are isomorphic as $\Mp$-module. Together with Proposition \ref{P3.13}, we
have following proposition.
\bp{P4.12} Let $M^{i}$ $(i=1,2,3)$ be irreducible
$\Mp$-modules. Then the fusion rule $N_{ M^{1}M^{2}}^{M^{3}}$ is
invariant under the permutation of $\{1,2,3\}$.
\ep

\subsection{Some useful lemmas}\label{S2.3}
\bl{L3.14}Let $(M,Y_{M})$ be an $\N$-gradable $V$-module. Then following hold
for any $a\in V$ homogeneous, $u\in M$ and homogeneous vector $v\in M$ of weight
$\wt(v)$:

\noindent
(1) For all $m,n\in\N$ such that $n\geq m$,
\eqa 
\Res _{z} {(1+z)^{\swt(a)+m}\over z^{2+n}}Y_{M}(a,z)u\in O(M).\nn
\eeqa

\noindent
(2) For any $n\in\Z_{>0}$, $ 
[L(-n)v]=(-1)^{n-1}[\ome*v-n v*\ome - \wt(v) v]$.
\el
\proof Part (1) is shown in Lemma 1.5.3 in \cite{FZ}. For (2), we can show that
$$[L(-n)v]=(-1)^{n-1}[(L(-1)-(n-1)(L(-2)+L(-1)))v]$$
for all $u\in M$ and $n\in\Z_{>0}$ by induction on $n$ and (1). Hence
\eqa[L(-n)v]&=&(-1)^{n-1}[(L(-2)+2L(-1)+L(0)-n(L(-2)+L(-1))-L(0))v]\nn\\
&=&(-1)^{n-1}[\ome*v-n v*\ome - \wt(v) v].\nn \eeqa\qed

\bl{L5.1} Let $M$ be a $V$-module and let $U$ be a subspace of $M$ such that
$a*U\subset U$ and $U*a\subset U$ hold for any $a\in V$. If $a\in V$ homogeneous
of positive weight and $u\in U$ satisfy that
$a_{n}u\in U+O(M)$ for any
$1\leq n\leq \wt(a)-1$, then $a_{n}u\in U+O(M)$ for all 
$n\leq\wt(a)-1$. In particular, if $u\in M$ is homogeneous, then
$L(-n)u\in U+O(M)$ for any $n\in\N$.
\el
\proof We show that 
\eqa a_{-n}u\in U+O(M), \label{5.1}\eeqa
for any $n\in\N$ by induction on $n$. From the assumption of this lemma, we have 
$$a_{0}u=a*u-u*a-\sum_{i=1}^{\swt(a)-1} {\wt(a)-1\choose i}
a_{i}u\in U+O(M).$$
So (\ref{5.1}) holds for $n=0$. In the case $n=1$, we have 
$$
a_{-1}u=u*a-\sum_{i=1}^{\swt(a)-1} {\wt(a)-1\choose i} a_{i-1}u.$$
Since $a_{n}u\in U+O(M)$ for every $0\leq n\leq\wt(a)-1$, we see that
(\ref{5.1}) holds for $n=1$. Assume that $\ell\in\Z_{>0}$ and (\ref{5.1}) holds for
any
$0\leq n\leq \ell$. Then by Lemma \ref{L3.14} (1), 
$$a_{-\ell-1}u+\sum_{i=1}^{\swt(a)-1} {\wt(a)\choose i}
a_{i-\ell-1}u\in O(M).$$
By induction hypothesis, we have $a_{-n-1}u\in U+O(M)$.\qed
\vs

Let $M$ be a $V$-module. We define a linear endomorphism
$\Ome_{M}:M\rightarrow M$ by
$$\Ome_{M}(u)=e^{L(1)}e^{\pi i L(0)} u$$
for $u\in M$. We remark 
that the operator $\Ome_{V}$ induces an anti-automorphism on $A(V)$ given in
\cite{Z} (see also \cite{DLM1}). A proof is similar to that of \cite{DLM1}.
\bp{P3.5}The linear map $\Ome_{M}:M\rightarrow M$ satisfies following
properties:
\eqa
\Ome_{M}(a*u)&=&\Ome_{M}(u)*\Ome_{V}(a),\
\Ome_{M}(u*a)=\Ome_{V}(a)*\Ome_{M}(u),\nn\\
\Ome_{M}(a\circ u)&=&-\Ome_{V}(a)\circ\Ome_{M}(u)\nn
\eeqa
for any $a\in V$ and $u\in M$.
\ep
\proof One can find in
\cite{FHL} the following conjugation formulas:
\eqa
z_{1}^{L(0)} Y_{M}(a,z_{0})z_{1}^{-L(0)}&=&Y_{M}(z_{1}^{L(0)}
a,z_{1}z_{0}),\nn\\
e^{z_{1}L(1)}Y_{M}(a,z_{0})e^{-z_{1}L(1)}&=&Y_{M}
\left(e^{z_{1}(1-z_{1}z_{0})L(1)}(1-z_{1}z_{0})^{-2L(0)}a,{z_{0}\over
1-z_{1}z_{0}}\right)\nn
\eeqa
for every $a\in V$. So we have  
\eqa 
&&\Res_{z_{0}}{(1+z_{0})^{\swt(a)+n}\over z_{0}^{m}}e^{L(1)}e^{\pi i
L(0)}Y_{M}(a,z_{0}) u\nn\\
&&{}=\Res_{z_{0}}{(1+z_{0})^{\swt(a)+n}\over
z_{0}^{m}}Y_{M}\left(e^{(1+z_{0})L(1)}(1+z_{0})^{-2\swt(a)}e^{\pi
i L(0)}a,{-z_{0}\over 1+z_{0}}\right)\Ome_{M}(u).\nn
\eeqa
for every $m,n\in \Z$, $a\in V$and $u\in M$. Here we replace
$-z_{0}/(1+z_{0})$ with
$w$ and apply the formula for change of variables (see \cite{Z}):
$$\Res_{w} g(w) = \Res_{z} (g(f(z)){d\over dz}f(z)),
$$
where $g(w)\in M((w))$ and $f(z)\in \C[[z]]$. Then we have 
\eqa
&&\Res_{z_{0}}{(1+z_{0})^{\swt(a)+n}\over z_{0}^{m}}e^{L(1)}e^{\pi i
L(0)}Y_{M}(a,z_{0}) u\nn\\
&&{}=(-1)^{m+1}\Res_{w}{(1+w)^{\swt(a) -n+m-2}\over w^{m}}
Y_{M}\left(\sum_{k=0}^{\infty}{(1+w)^{-k}\over k!}L(1)^{k}
e^{\pi i L(0)}a,w\right)\Ome_{M}(u)\nn\\
&&{}=(-1)^{m+1}\sum_{k=0}^{\infty}\Res_{w}{(1+w)^{\swt(a)
-n+m-2-k}\over w^{m}} Y_{M}\left({L(1)^{k}\over k!}
e^{\pi i L(0)}a,w\right)\Ome_{M}(u).\nn
\eeqa
Hence if we take  $m=1$ and $n=0$, then
 \eqa\Ome_{M}(a*u)&=&\sum_{k=0}^{\infty}\Res_{w}{(1+w)^{\swt(a) -1-k}\over w}
Y_{M}\left({L(1)^{k}\over k!} e^{\pi i L(0)}a,w\right)\Ome_{M}(u)\nn\\
&=&\sum_{k=0}^{\infty}\Ome_{M}(u)*\left({L(1)^{k}\over k!}
e^{\pi i L(0)}a\right)\nn\\
&=&\Ome_{M}(u)*\Ome_{V}(a).\nn
\eeqa
Similarly if we take $m=1$ and $n=-1$, then
$\Ome_{M}(u*a)=\Ome_{V}(a)*\Ome_{M}(u),$ and if $m=2$ and $k=0$, then
$\Ome_{M}(a\circ u)=-\Ome_{V}(a)\circ \Ome_{M}(u).$
\qed
\vs

\noindent
Thus $\Ome_{M}$ induces a linear
map on $A(M)$ (also denoted by $\Ome_{M}$) such that
$\Ome_{M}([a*u])=\Ome_{M}([u])*\Ome_{V}([a]),\
\Ome_{M}([u*a])=\Ome_{V}([a])*\Ome_{M}([u])$ 
hold for all $a\in V$ and $u\in M$. 

Let $M$ be a weak $V$-module. An element $u\in M$ is called a {\it
lowest weight vector of weight $\wt(u)\in\C$} if $L(n)u=\wt(u)\delta_{n,0} u$ for
any $n\in\N$. If $a\in V$ is a lowest weight
vector, we have the following commutation relation:
\eqa
[L(m),a_{n}]=((\wt(a)-1)(m+1)-n)a_{m+n},\ {\rm for}\ m,n\in\Z.\label{2.27}
\eeqa
\bl{L2.8} Let $M$ be a $V$-module, and let $a\in V$ and
$v\in M$ be lowest weight vectors. Then the subset of M spanned by vectors 
$a_{n_{1}}a_{n_{2}}\cdots a_{n_{k}}v,$ $k\in\Z_{>0}$ and
$n_{i}\in\Z$ $( i=1,2,\ldots,k)$ is invariant under the action of $L(n)$ for
$n\in\N$.
\el
\proof Set $U=\Sp\{a_{n_{1}}a_{n_{2}}\cdots a_{n_{k}}v\ |\
k\in\Z_{>0},\ n_{i}\in\Z\ \}$. To show $L(n)U\subset U$ for all
$n\in\N$, we prove that $
L(n)a_{n_{1}}a_{n_{2}}\cdots a_{n_{k}}v\in U$ 
for any $n_{i}\in\Z$. By the commutation relation (\ref{2.27}), we have 
\eqa
L(n)a_{n_{1}}a_{n_{2}}\cdots a_{n_{k}}v&=&\sum_{i=1}^{k}
a_{n_{1}}\cdots [L(n),a_{n_{i}}]\cdots
a_{n_{k}}v+a_{n_{1}}a_{n_{2}}\cdots a_{n_{k}}L(n)v\nn\\
&=&\sum_{i=1}^{k}((\wt(a)-1)(n+1)-n_{i})
a_{n_{1}}\cdots a_{n_{i+n}}\cdots
a_{n_{k}}v\nn\\
&&{}+a_{n_{1}}a_{n_{2}}\cdots a_{n_{k}}L(n)v.\label{2.29}
\eeqa
Since $L(n)v=\wt (v)\delta_{n,0}\, v$ for $n\in\N$, the
right-hand side of (\ref{2.29}) is in $U$.\qed
\bl{L2.9} Let $M$, $a$ and $v$ be as in Lemma $\ref{L2.8}$.
Then the subspace of M spanned by vectors
\eq
L(-m_{1})L(-m_{2})\cdots L(-m_{k})a_{n_{1}}a_{n_{2}}\cdots
a_{n_{\ell}}v,\ k,\ \ell\in\N,\ m_{i}\in\Z_{>0},\
n_{j}\in\Z,\label{2.30}
\eeq
is invariant under the actions of $L(n)$ $(n\in\Z)$ and $a_{m}$
$(m\in\Z)$.
\el
\proof Set $U$ be the subspace of $M$ spanned by vectors (\ref{2.30}). By
Poincar{\'e}-Birkhoff-Witt (PBW) Theorem for Virasoro algebra and Lemma \ref{L2.8},
we see that $U$ is invariant under the action of $L(n)$ for $n\in\Z$.

Using induction on k, we prove that 
\eq
a_{m}L(-m_{1})L(-m_{2})\cdots L(-m_{k})a_{n_{1}}a_{n_{2}}\cdots
a_{n_{\ell}}v\in U,\label{2.31}
\eeq
for any $k,\ell\in\N,$ $m,m_{i}\in\Z_{>0}$ and
$n_{j}\in\Z$. The case $k=0$ is clear. Assume that (\ref{2.31})
holds for $k=p\in\N$. Then if we put 
$v'=L(-m_{2})\cdots L(-m_{p+1}) a_{n_{1}}a_{n_{2}}\cdots
a_{n_{\ell}}v,$
we have
\eqa 
a_{m}L(-m_{1})v'&=&-[L(-m_{1}), a_{m}]v'+L(-m_{1})a_{m} v'\nn\\
&=&((\wt(a)-1)(m_{1}+1)+m)a_{n-m_{1}}v'+L(-m_{1})a_{m}
v'.\nn
\eeqa
By induction hypothesis and the previous paragraph, $a_{n-m_{1}}v'$, $a_{m} v'$ and
$L(-m_{1})a_{m} v'$ belong to $U$. Hence (\ref{2.31}) holds for $k=p+1$.\qed

\subsection{A spanning set of $A(M)$}\label{S4.2}

Let $M$ be an irreducible $\Mp$-module. Then from Lemma \ref{L2.9}, the subspace
of $M$ spanned by vectors
\eqa L(-m_{1})L(-m_{2})\cdots L(-m_{k})J_{n_{1}}J_{n_{2}}\cdots
J_{n_{\ell}}v_{M}\label{5.2}\eeqa
for $k, \ell\in \N$, $m_{i}\in\Z_{>0}$ and $n_{j}\in \Z$ is
invariant under the actions of $L(m)$ and $J_{n}$ for all $m,n\in\Z$, where $v_{M}$
is the lowest weight vector of $M$ given in Table \ref{Table4.1}. Since $\Mp$ is
generated by
$\ome$ and $J$, this subspace is invariant under the action of
$\Mp$. Thus this is $\Mp$-submodule of $M$. By the irreducibility of
$M$, we have following lemma.
\bl{L5.2} An irreducible $\Mp$-module $M$ is spanned by vectors $(\ref{5.2})$.
\el
One of main results of this subsection is the following proposition. 
\bp{P5.3} Let $M$ be an irreducible $\Mp$-module. Suppose that
there exists a lowest weight vector $u\in M$ such that for any
$n\in\Z_{>0}$, $J_{n}v_{M}$ and $J_{n}u$ belong to the submodule for
Virasoro algebra of $M$ generated by $v_{M}$ and $u$. Then $M$
is spanned by the set 
\eqa \{a*v_{M}*b, a*u*b| a,b\in\Mp\}+O(M).\label{5.3}\eeqa 
\ep
To prove the proposition, we need some lemmas. The proof of proposition is given
after Lemma \ref{L5.5}.
\bl{L4.6} Let $M$ be an $\Mp$-module. Then for any
$m,n\in\Z$, the commutator $[J_{m},J_{n}]$ is a linear combination of 
\eqa L(m_{1}) L(m_{2})\cdots L(m_{p})J_{m_{p+1}}\  {\it and}\ 
L(n_{1})L(n_{2})\cdots L(n_{q}),\ {\it for}\ m,n\in\Z.\nn\eeqa
\el
\proof This is proved by the similar argument of Lemma of \cite{DN}.\qed

\bl{L5.4}Let $u\in M$ be a lowest weight vector, and set $U_{M,u}$ to be the
subspace of $M$ spanned by the set $(\ref{5.3})$. Fix a positive integer
$\ell$. Suppose that for any $1\leq k \leq \ell$ and
$n_{1},\ldots,n_{k}\in\Z$, $J_{n_{1}}\cdots J_{n_{k}}v_{M}$ and
$J_{n_{1}}\cdots J_{n_{k}}u$ lie in $U_{M,u}$. Then for any $1\leq
k\leq \ell$, $1\leq j\leq k$, $s\in\Z_{>0}$ and $p_{1},\ldots,
p_{s}\in\Z$, we have 
$$ J_{n_{1}}\cdots J_{n_{j}}L(p_{1})\cdots
L(p_{s})J_{n_{j+1}}\cdots J_{n_{k}}v_{M},
J_{n_{1}}\cdots J_{n_{j}}L(p_{1})\cdots
L(p_{s})J_{n_{j+1}}\cdots J_{n_{k}}u\in U_{M,u}.$$ 
\el
\proof Let $v$ be either $v_{M}$ or $u$. We have to show that
\eqa J_{n_{1}}\cdots J_{n_{j-1}}L(p_{1})\cdots
L(p_{s})J_{n_{j}}\cdots J_{n_{k}}v\in U_{M,u},\label{5.6}\eeqa
for any $s\in\Z_{>0}$ and $n_{1},\ldots,n_{k},p_{1},\ldots,p_{s}\in\Z$. By
PBW Theorem for Virasoro algebra, we may assume that
$p_{1}\leq p_{2}\leq\cdots\leq p_{s}$. If $p_{1}\geq 0$, then
$p_{s}\geq0$. Then we view (\ref{5.6}) as follows:
\eqa&& J_{n_{1}}\cdots J_{n_{j-1}}L(p_{1})\cdots
L(p_{s})J_{n_{j}}\cdots J_{n_{k}}v\nn\\
&&\ \ {}=\sum_{i=j}^{k}
J_{n_{1}}\cdots J_{n_{j-1}}L(p_{1})\cdots
L(p_{s-1})J_{n_{j}}\cdots[L(p_{s}),J_{n_{i}}]\cdots 
J_{n_{k}}v\nn\\
&&\ \ \ \ \ \ {}+ J_{n_{1}}\cdots J_{n_{j-1}}L(p_{1})\cdots
L(p_{s-1})J_{n_{j}}\cdots J_{n_{k}}L(p_{s}) v\nn\\
&&\ \ {}=\sum_{i=j}^{k} (3(p_{s}+1)-n_{i})J_{n_{1}}\cdots
J_{n_{j-1}}L(p_{1})\cdots L(p_{s-1})J_{n_{j}}\cdots 
J_{p_{s}+n_{i}}\cdots  J_{n_{k}}v\nn\\
&&\ \ \ \ \ \ {}+\wt(v)\del_{p_{s},0} J_{n_{1}}\cdots
J_{n_{j-1}}L(p_{1})\cdots L(p_{s-1})\cdots J_{n_{k}}v.\nn
\eeqa
If $p_{1}<0$, then we view as follows:  
\eqa &&J_{n_{1}}\cdots J_{n_{j-1}}L(p_{1})\cdots
L(p_{s})J_{n_{j}}\cdots J_{n_{k}}v\nn\\
&&\ \ {}=-\sum_{i=1}^{j-1}
J_{n_{1}}\cdots[L(p_{1}),J_{n_{i}}]\cdots
J_{n_{j-1}}L(p_{2})\cdots L(p_{s})J_{n_{j}}\cdots  J_{n_{k}}v\nn\\
&&\ \ \ \ \ \ {}+ L(p_{1})J_{n_{1}}\cdots
J_{n_{j-1}}L(p_{2})\cdots L(p_{s})J_{n_{j}}\cdots  J_{n_{k}}v.\nn\\
&&\ \ {}=-\sum_{i=1}^{j-1} (3(p_{1}+1)-n_{i})J_{n_{1}}\cdots J_{p_{1}+n_{i}}\cdots 
J_{n_{j-1}}L(p_{2})\cdots L(p_{s})J_{n_{j}}\cdots  J_{n_{k}}v\nn\\
&&\ \ \ \ \ \ {}+ L(p_{1})J_{n_{1}}\cdots
J_{n_{j-1}}L(p_{2})\cdots L(p_{s})J_{n_{j}}\cdots  J_{n_{k}}v.\nn
\eeqa 
Hence by using induction on $s$,  we see that (\ref{5.6}) follows from Lemma
\ref{L5.1}.
\qed
\bl{L5.5} Let $M$ and $u$ be as in Proposition $\ref{P5.3}$.
Then for any positive integer $\ell$ and $n_{i}\in\Z$ $(i=1,\ldots,\ell)$,
we have $J_{n_{1}}\cdots J_{n_{\ell}}v_{M},J_{n_{1}}\cdots J_{n_{\ell}}u\in
U_{M,u}$, where $U_{M,u}$ is the subset of $M$ defined in Lemma
$\ref{L5.4}$.
\el
\proof Let $v$ be either $v_{M}$ or $u$.
By using induction on $\ell$, we prove
\eqa J_{n_{1}}\cdots J_{n_{\ell}}v\in U_{M,u},\label{5.7}\eeqa
 for any $n_{i}\in\Z$. In the case $\ell=1$, the assumptions of this lemma
implies that for any positive integer $n$, $J_{n}v$ is a linear combination of
$L(-m'_{1})\cdots L(-m'_{p})v_{M}$ and
$L(-n'_{1})\cdots L(-n'_{q})u$ for $m'_{i},n'_{j}\in\Z_{>0}$. Hence by
Lemma \ref{L3.14} (2), $J_{n}u\in U_{M,u}$ holds for any $n\in\Z_{>0}$.
Thus (\ref{5.7}) for $\ell=1$ follows from Lemma \ref{L5.1}. Assume
that $p\in\Z_{>0}$ and that (\ref{5.7}) holds for any $\ell\leq p$. Then we
have to show (\ref{5.7}) for $\ell=p+1$. By induction hypothesis, $J_{n_{2}}\cdots
J_{n_{p+1}}v\in U_{M,u}$. Hence by Lemma
\ref{L5.1}, it is sufficient to prove that  (\ref{5.7}) for $\ell=p+1$ holds for any
$n_{1}\in\Z_{>0}$. Then we have 
\eqa J_{n_{1}}\cdots J_{n_{p+1}}v=\sum_{i=2}^{p+1}
J_{n_{2}}\cdots [J_{n_{1}}, J_{n_{i}}]\cdots J_{n_{p+1}}v+J_{n_{2}}\cdots
J_{n_{p+1}}J_{n_{1}}v.\label{5.9}\eeqa By Lemma \ref{L4.6}, $[J_{n_{1}}, J_{n_{i}}]$ is
expressed by a linear combination of $L(m'_{1})\cdots L(m'_{q})J_{m'_{q+1}}$ and
$L(n'_{1})\cdots L(n'_{r})$ for $m'_{i},n'_{j}\in \Z$, and
 $J_{n_{1}}v$ can be expressed by a
linear combination of 
$L(-m''_{1})\cdots L(-m''_{s})v_{M}$ and
$L(-n''_{1})\cdots L(-n''_{t})u$ for $m''_{i},n''_{j}\in\Z_{>0}$ by assumption of this
Lemma. Hence by induction hypothesis and Lemma \ref{L5.4} , the right-hand
side of the equality (\ref{5.9}) lies in $U_{M,u}$.\qed
\vs

Now we prove Proposition \ref{P5.3}.

\noindent
{\it Proof of Proposition \ref{P5.3}.} Let $U_{M,u}$ be the same set as in Lemma
\ref{L5.4}. Then we show $M=U_{M,u}$. By Lemma
\ref{L5.2}, it is enough to show that 
\eqa L(-m_{1})L(-m_{2})\cdots L(-m_{k})J_{n_{1}}J_{n_{2}}\cdots
J_{n_{\ell}}v_{M}\in U_{M,u}\label{5.10}\eeqa
for any $m_{i}\in\Z_{>0}$ and $n_{j}\in\Z$. But Lemma \ref{L5.5}
tells us that $J_{n_{1}}\cdots J_{n_{\ell}}v_{M}\in U_{M,u}$, and since
the vectors of this form is homogeneous,  (\ref{5.10}) follows from
Lemma \ref{L5.1}.\qed
\vs

We consider the case $M=\M^{\pm}$, $\Ml$ ($\lm\neq0,1/2$) and
$\Mtp$. Since $\deg(J_{n}v_{M})=3-n\leq 2$ if $n\geq 1$, by the
irreducible decomposition (\ref{4.20}), (\ref{4.22}) and
(\ref{4.31}), $J_{n}v_{M}$ ($n\geq1$) lies in the submodule for Virasoro algebra of
$M$ generated by $v_{M}$. Thus if we take
$u=v_{M}$ in Proposition \ref{P5.3}, we have  
$M=\Sp\, \{a*v_{M}*b\,|\,a,b\in\Mp\}+O(M),$ then $A(M)$ is generated by $[v_{M}]$ as
$A(\Mp)$-bimodule.

In the case $M=\Ml$ ($\lm^{2}=1/2$), by (\ref{4.20}) we have the irreducible
decomposition $ M=L\left(1,{1/4}\right)\ops L\left(1,{9/4}\right)\Ops
L\left(1,{25/4}\right)\ops\cdots.$ Now let $u$ be the lowest weight vector of
$L(1,9/4)$, then $\deg (J_{n}v_{M})=3-n\leq 2$ and
$\deg(J_{n}u)=5-n\leq 4$ for any positive integer $n$. Hence $J_{n}v_{M}$ and
$J_{n}u$ ($n\geq1$) lie in $L(1,1/4)\ops L(1,9/4)$. Thus Proposition
\ref{P5.3} implies that $\Ml=\Sp \{a*v_{M}*b,a*u*b|a,b\in\Mp\}+O(M)$. In this case,
$A(\Ml)$ is generated by $[v_{M}]$ and $[u]$ as
$A(\Mp)$-bimodule. 

In the remaining case $M=\Mtm$, by (\ref{4.32}), we have the irreducible
decomposition $M=L\left(1,{9/16}\right)\ops L\left(1,{25/ 16}\right)\ops
 L\left(1,{121/16}\right)\ops\cdots.$ Let $u'$ be the lowest weight
vector of $L(1,25/16)$, then $\deg(J_{n}v_{M})=3-n\leq 2$ and
$\deg(J_{n}u)=4-n\leq 3$ for any positive integer $n$. Hence $J_{n}v_{M}$ and
$J_{n}u'$ ($n\geq1$) are in $L(1,9/16)\ops L(1,25/16)$. Thus 
by Proposition
\ref{P5.3}, $\Mtm=\Sp \{a*v_{M}*b,a*u'*b|a,b\in\Mp\}+O(M)$, and $A(\Mtm)$ is
generated by $[v_{M}]$ and $[u']$ as $A(\Mp)$-bimodule. Summarizing, we have the
following proposition.
\bp{P5.6} Let $M$ be a irreducible $\Mp$-module. Then following
hold.

\noindent
(1) If  $M=\M^{\pm}$, $\Ml$ $(\lm^{2}\neq0,1/2)$ or $\Mtp$, then
$A(M)$ is generated by $[v_{M}]$ as $A(\Mp)$-bimodule.

\noindent
(2) If $M=\Ml$ $(\lm^{2}=1/2)$, then $A(M)$ is generated by $[v_{M}]$
and $[u]$ as $A(\Mp)$-bimodule, where $u$ is a lowest weight
vector of weight $9/4$.

\noindent
(3) If $M=\Mtm$, then then $A(M)$ is generated by $[v_{M}]$
and $[u']$ as $A(\Mp)$-bimodule, where $u'$ is a lowest weight
vector of weight $25/16$.
\ep
\br{R5.7} If $M$ is the vertex operator algebra $\Mp$ itself, this
result is clear because one have $a*\1=a=\1*a$ for all $a\in\Mp$.
\er

As a corollary of Proposition \ref{P5.6}, we have:
\bc{C5.8}Let $M$, $N$, $L$ be irreducible $\Mp$-modules.

\noindent
(1)  If  $M=\M^{\pm}$, $\Ml$ $(\lm^{2}\neq0,1/2)$ or $\Mtp$, then 
$N_{MN}^{L}\leq 1.$

\noindent
(2) If $M=\Ml$ $(\lm^{2}=1/2)$ or $\Mtm$, then  $N_{MN}^{L}\leq 2.$
\ec
\proof If $A(M)$ is generated by $n$ vectors
$[v_{1}]$,$[v_{2}],\ldots,[v_{n}]$ as $A(\Mp)$-bimodule, then the
contraction $(L_{0})^{*}\cdotp A(M)\cdotp N_{0}$ is spanned by $n$
vectors $ v_{L}'\ots [v_{1}]\ots v_{N},\ldots,v_{L}'\ots [v_{n}]\ots
v_{N}.$ 
Hence we have $\dim_{{\sC}}(L_{0})^{*}\cdotp A(M)\cdotp
N_{0}\leq n$. On the other hand, by Proposition \ref{P3.12},
\eqa N_{MN}^{L}\leq \dim_{{\sC}}((L_{0})^{*}\cdotp A(M)\cdotp
N_{0})^{*}= \dim_{{\sC}}(L_{0})^{*}\cdotp A(M)\cdotp
N_{0}\leq n.\nn\eeqa
Thus the corollary follows from Proposition \ref{P5.6}.\qed


\section{Fusion rules}\label{S5}
This section is consisted of two subsections. In first subsection we
prove some propositions which give triples of irreducible $\Mp$-modules whose fusion
rules are nonzero, and state main theorem. In second subsection, we
prove the main theorem.
\subsection{Main theorem}\label{S5.1}
Recall the linear map $ I_{\lm\mu}:\Ml\to z^{\lm\mu}
\Hom(\Mu,M(1,\lm+\mu))[[z,z^{-1}]]$ defined in (\ref{4.3}) and (\ref{4.4})
for $\lm,\mu\in\C$. By the arguments of Section 8.6 of \cite{FLM}, we see that
$I_{\lm\mu}$ is an intertwining operator of type
\fusion{\Ml}{\Mu}{ M(1,\lm+\mu)}. Since $e^{\mu}\in\Mu_{0}$, we have 
\eqa I_{\lm\mu}(e^{\lm},z)e^{\mu}=z^{\lm\mu} \exp\left(\sum_{n>0}
{\lm h(-n)\over n} z^{n}\right)e^{\lm+\mu},\nn\eeqa
and its coefficient of $z^{\lm\nu}$ is $e^{\lm+\mu}$ which is
nonzero. Thus this intertwining operator $I_{\lm\mu}$ is nonzero.
This implies that the fusion rule
$N_{\ss{ \Ml\Mu}}^{M(1,\lm+\mu)}$ is nonzero for any
$\lm,\mu\in\C$. By Proposition \ref{P4.3}, there exists an isomorphism
$f$ from $\Mu$ onto $M(1,-\mu)$ of $\Mp$-modules. We define a linear map
\eqa Y_{\lm}\circ f:\Ml&\to&
z^{-\lm\mu}(\Hom(\Mu,M(1,\lm-\mu)))[[z,z^{-1}]],\nn\\
u\ \ \ &\mapsto&Y_{\lm}\circ f (u,z)=Y_{\lm}(u,z)\circ f.\nn\eeqa
Then we can easily see that $Y_{\lm}\circ f$ gives a nonzero
intertwining operator of type  \fusion{\Ml}{\Mu}{ M(1,\lm-\mu)}. Thus following
proposition holds.
\bp{P6.1}For $\lm,\mu,\nu\in\C$, the fusion rule $N_{\ss{
\Ml\Mu}}^{\ss{ \Mn}}$ is nonzero if $\nu^{2}=(\lm\pm\mu)^{2}$.
\ep
In particular if $\lm=0$, we have that $N_{\ss{ \M\Mu}}^{\ss{\Mn}}$
is nonzero if $\mu^{2}=\nu^{2}$. In fact, we have the following proposition.
\bp{P6.2} For $\mu,\nu\in\C$ such that $\mu^{2}=\nu^{2}$ , the
fusion rules $N_{\ss{ \M^{\pm}\Mu}}^{\ss{\Mn}}$ are nonzero.
Furthermore the fusion rules 
$N_{\ss{\Mp\M^{\pm}}}^{\ss{\M^{\pm}}}$ and  $N_{\ss{
\Mm\M^{\pm}}}^{\ss{\M^{\mp}}}$ are nonzero.
\ep
\proof Since $(\Mu,I_{0\mu})$ is an irreducible $\M$-module, $I_{0\mu}$ gives
nonzero intertwining operators of types \fusion{\M^{\pm}}{\Mu}{\Mu} by Proposition
11.9 of \cite{DL}. Because
$M(1,-\mu)$ is isomorphic to $\Mu$ as $\Mp$-module, the first statement
holds. If $\mu=0$, then $Y=I_{00}$ gives nonzero intertwining operators of types  
\fusion{\M^{\pm}}{\M}{\M}. Since $\th$ is an automorphism of $\M$, we
have $\th Y(a,z)\th =Y(\th(a),z)$  for every $a\in\M$. This implies
that $ Y(a,z)\M^{\pm}\subset \M^{\pm}((z))$ for every $a\in\Mp$ and
$Y(a,z)\M^{\pm}\subset \M^{\mp}((z))$ for every $a\in\Mm$. Thus $Y$
gives nonzero intertwining operators of type 
\fusion{\Mp}{\M^{\pm}}{\M^{\pm}} and \fusion{\Mm}{\M^{\pm}}{\M^{\mp}}.\qed
\vs

Next we recall the linear map $ I_{\lm}^{\th}:\Ml\to (\End \Mt)\{z\}$ 
defined in (\ref{4.125}) and (\ref{4.13}). The arguments of Chapter 9 of \cite{FLM}
shows that $I_{\lm}^{\th}$ is an intertwining operator of type 
\fusion{\Ml}{\Mt}{\Mt}. Let $p_{\pm}$ be projections from $\Mt$ onto $\Mt^{\pm}$
and let $\iota_{\pm}$ be inclusions from $\Mt^{\pm}$ into $\Mt$ respectively. Define
a linear map
$p_{\beta}\circ I_{\lm}^{\th}\circ \iota_{\al}$ by
\eqa 
p_{\beta}\circ I_{\lm}^{\th}\circ \iota_{\al}:\Ml&\to&
(\Hom(\Mt^{\al},\Mt^{\beta}))\{z\},\nn\\
u&\mapsto&p_{\beta}\circ
I_{\lm}^{\th}\circ \iota_{\al}(u,z)=p_{\beta}(I_{\lm}^{\th}(u,z)
 \iota_{\al}),\nn\eeqa where $\al,\beta$ are $+$ or $-$. Then we
have following lemma.
\bl{L6.3} For any $\al,\beta\in\{+,-\}$, $p_{\beta}\circ I_{\lm}^{\th}\circ
\iota_{\al}$ is a nonzero intertwining operator of type \fusion{\Ml}{\Mt^{\al}}{
\Mt^{\beta}} if $\lm\neq 0$.
\el
\proof It is clear that $p_{\beta}\circ I_{\lm}^{\th}\circ \iota_{\al}$ is intertwining
operator of type  \fusion{\Ml}{\Mt^{\al}}{ \Mt^{\beta}}. We next show
$p_{\beta}\circ I_{\lm}^{\th}\circ \iota_{\al}$ is nonzero if
$\lm\neq 0$. By direct calculation, we can see that the coefficient of $z^{-\lm^{2}/2}$
in $I_{\lm}^{\th}(e^{\lm},z)1$ is $1$ and that of $z^{-\lm^{2}/2+1/2}$ is
$\lm h(-1/2)1$, and the coefficient of $z^{-\lm^{2}/2-1/2}$ in
$I_{\lm}^{\th}(e^{\lm},z)h(-1/2)1$  is
$-\lm 1$ and that of $z^{-\lm^{2}/2+1}$ is
$-(2/3)\lm^{2}h(-{3/2})1+2 \lm^{2}(1-{3/2}\lm^{2}
)h(-{1/2})^{3}1$. If $\lm\neq 0$, all of these coefficients
are nonzero. This proves that $p_{\beta}\circ I_{\lm}^{\th}\circ
\iota_{\al}(e^{\lm},z)$ is nonzero for all $\al,\beta\in\{+,-\}$ if $\lm$ is
nonzero.\qed
\vs

\bp{P6.4}(1)  For any nonzero $\lm\in\C$ and
$\al,\beta\in\{+,-\}$, $N_{\ss
{\Ml\Mt^{\al}}}^{\ss{\Mt^{\beta}}}$ is nonzero.

\noindent
(2) The fusion rules $N_{\ss{
\Mp\Mt^{\pm}}}^{\ss{\Mt^{\pm}}}$ and $N_{\ss{
\Mm\Mt^{\pm}}}^{\ss{\Mt^{\mp}}}$ are nonzero.
\ep
\proof The part (1) follows from Lemma \ref{L6.3}. By the definition of the
action of $\th$ on $\Mt$ (\ref{4.16}),  we can show that $\th
I_{0}^{\th}(a,z)\th=I_{0}^{\th}(\th(a),z)$ holds for every 
$a\in\M$. This implies that 
\eqa I_{0}^{\th}(a,z)\Mt^{\pm}&\subset&\Mt^{\pm}((z))\ \ 
{\rm if}\ a\in\Mp,\nn\\
 I_{0}^{\th}(a,z)\Mt^{\pm}&\subset&\Mt^{\mp}((z))\ \ 
{\rm if}\ a\in\Mm.\nn\eeqa
Because $I_{0}^{\th}(\1,z)=\id_{\ss{ \Mt}}$ and
$I_{0}^{\th}(h(-1)\1,z)=h(z)$, we see that $p_{\pm}\circ I_{0}^{\th}\circ
\iota_{\pm}$ ($p_{\pm}\circ I_{0}^{\th}\circ
\iota_{\mp}$) give nonzero intertwining operators of type
\fusion{\Mp}{\Mt^{\pm}}{\Mt^{\pm}} (resp. \fusion{\Mm}{\Mt^{\pm}}{\Mt^{\mp}})
. This shows the part (2).\qed
\vs

The series of Propositions \ref{P6.1}, \ref{P6.2} and \ref{P6.4}
together with Proposition \ref{P4.12} give us triples $(M,N,L)$ of irreducible
$\Mp$-modules which the fusion rule $N_{MN}^{L}$ is nonzero. In fact, we
have following theorem.
\bt{T6.5} $({\rm Main Theorem})$ Let $M$, $N$ and $L$ be irreducible
$\Mp$-modules.

\noindent
(1) If $M=\Mp$, then $N_{\ss{\Mp}{ N}}^{L}=\delta_{N,L}$.

\noindent
(2) If $M=\Mm$, then $N_{\ss{\Mm}{ N}}^{L}$ is $0$ or $1$, and
$N_{\ss{\Mm}{ N}}^{L}=1$ if and only if the pair $(N,L)$ is one of
following pairs:
\eqa(\M^{\pm},\M^{\mp}),\ (\Mt^{\pm},\Mt^{\mp}),\ (\Ml,\Mu)\
(\lm^{2}=\mu^{2}).\nn
\eeqa 

\noindent
(3) If $M=\Ml$ ($\lm\neq0$), then $N_{\ss{\Ml}{ N}}^{L}$ is $0$ or $1$,  and
$N_{\ss{\Ml }{N}}^{L}=1$ if and only if the pair $(N,L)$ is one of
following pairs:
\eqa&&(\M^{\pm},\Mu)\ (\lm^{2}=\mu^{2}),\ (\Mu,\Mn)\
(\nu^{2}=(\lm\pm\mu)^{2}),\nn\\
&&(\Mt^{\pm},\Mt^{\pm}),\
(\Mt^{\pm},\Mt^{\mp}).\nn
\eeqa

\noindent
(4) If $M=\Mtp$, then $N_{\ss{\Mtp}{ N}}^{L}$ is $0$ or $1$, and
$N_{\ss{\Mtp}{ N}}^{L}=1$ if and only if the pair $(N,L)$ is one of
following pairs:
\eqa&&(\M^{\pm},\Mt^{\pm}),\ (\Ml,\Mt^{\pm}).\nn
\eeqa

\noindent
(5) If $M=\Mtm$, then $N_{\ss{\Mtm} N}^{L}$ is $0$ or $1$, and
$N_{\ss{\Mtm} N}^{L}=1$ if and only if the pair $(N,L)$ is one of
following pairs:
\eqa&&(\M^{\pm},\Mt^{\mp}),\ (\Ml,\Mt^{\pm}).\nn
\eeqa

\et 

\subsection{Proof of main theorem}\label{S5.2}
Recall the direct sum decomposition $\Mp=L(1,0)\ops L(1,4)\ops \cdots$ (see
(\ref{4.22})). The lowest weight vector of $L(1,4)$ is $J$. Hence $h(-3)h(-1)\1$ is a
linear combination of
$L(-4)\1$, $L(-2)^{2}\1$ and $J$ (note that $L(-1)\1=0$). In fact,
$h(-3)h(-1)\1=(-1/9)(J-3L(-4)\1-4L(-2)^{2}\1)$. By Lemma \ref{L3.14} (2), we have
the following equality in $A(\Mp)$:
\eqa[J-4\ome^{*2}-17 \ome+9 h(-3)h(-1)\1]=0,\label{6.0}
\eeqa
where we use the notation
$[a]^{*n}=[a^{*n}]=[\overbrace{a*a*\cdots *a}^{n}]$ for
$a\in\Mp$ and $n\in\N$.

Now we prove the main theorem. We divide the proof into five steps Step 1-Step 5
where (1)-(5) proved respectively.

\noindent
{\it Proof of Theorem \ref{T6.5}.} By Theorem \ref{T4.8} and
Table \ref{Table4.1}, we see that two irreducible
$\Mp$-modules $N$ and $L$ are isomorphic each other if and
only if $a_{N}=a_{L}$ and $b_{N}=b_{L}$.

{\it Step 1.} If $M=\Mp$, then $v_{M}=\1$ (see Table
\ref{Table4.1}), and $(L_{0})^{*}\cdotp A(M)\cdotp
N_{0}=\C\ v_{L}'\ots [\1]\ots v_{N}$ by Proposition \ref{P5.6} (1). 
If the fusion rule $N_{MN}^{L}$ is nonzero, then $v_{L}'\ots
[\1]\ots v_{N}$ is also nonzero. On one hand, we have following
equalities:
\eqa a_{L}v_{L}'\ots [\1]\ots v_{N}&=&v_{L}'\ots [\ome]\ots
v_{N}= a_{N}v_{L}'\ots [\1]\ots v_{N},\nn\\
b_{L}v_{L}'\ots [\1]\ots v_{N}&=&v_{L}'\ots [J]\ots
v_{N}=b_{N}v_{L}'\ots [\1]\ots v_{N}.\nn\eeqa
This implies that if the fusion rule $N_{MN}^{L}$ is nonzero, then 
$a_{N}=a_{L}$ and $b_{N}=b_{L}$, hence $M$ and $N$ are
equivalent. Then (1) follows from Proposition \ref{P6.2} and Proposition 
\ref{P6.4}.

{\it Step 2.} Let $M=\Mm$. Then we have the irreducible
decomposition (\ref{4.22}). Hence $h(-3)h(-1)\1*v_{M}$ is in $L(1,1)$ and
it is linear combination of $L(-m_{1})\cdots L(-m_{k})v_{M}$ ($m_{i}\in
\Z_{>0}$, $m_{1}+\cdots+m_{k}\leq4$). In fact,
\eqa&&h(-3)h(-1)\1*v_{M}=3v_{M}+ 12L(-1)v_{M}+12L(-1)^{2}v_{M}
-8L(-3)v_{M}\ \ \ \ \ \ \ \ \ \ \ \ \ \ \nn\\
&&+16L(-2)L(-1)v_{M}-{1\over2}L(-4)v_{M}+{1\over4}L(-3)L(-1)v_{M}+{3\over2}
L(-2)L(-1)^{2}v_{M}.\label{6.3}\eeqa
 On the other hand, we have 
\eqa v_{L}'\ots [L(-m_{1})\cdots L(-m_{k})v_{M}]\ots v_{N}
=F(a_{L},a_{N}) v_{L}'\ots [v_{M}]\ots v_{N}\label{6.4}\eeqa 
for any $m_{i}\in\Z_{>0}$ by Lemma \ref{3.6}, where
$F\in\C[x,y]$ is given by
\eqa F=\prod_{i=1}^{k} (-1)^{m_{i}-1}(x- m_{i} y-\sum_{j=i+1}^{k}
m_{j}-\wt(v_{M})).\label{6.5}
\eeqa
Hence by (\ref{6.0}) and (\ref{6.3})-(\ref{6.5}), we have
\eqa &&v_{L}'\ots [(J-4\ome^{*2}-17 \ome+9
h(-3)h(-1)\1)*v_{M}]\ots v_{N} \nn\\
&&\ \ =f(a_{L},a_{N},b_{L}) v_{L}'\ots [v_{M}]\ots
v_{N}=0,\nn
\eeqa
where $f\in\C[x,y,z]$ is given by
$$f=z-4 x^{2}+x +{9\over 4} (x-y)(6 x^{2}-18 x y-12
y^{2}-21 x-23 y+11).$$
Proposition \ref{P3.5} and (\ref{6.0}) show that 
\eqa &&v_{L}'\ots [\Ome_{M}(v_{M})*\Ome_{\ss{\Mp}}(J-4\ome^{*2}-17 \ome+9
h(-3)h(-1)\1)]\ots v_{N} \nn\\
&&\ \ =v_{L}'\ots\Ome_{M}( [(J-4\ome^{*2}-17 \ome+9
h(-3)h(-1)\1)*v_{M}])\ots v_{N}\nn \\
&&\ \ =-f(a_{N},a_{L},b_{N}) v_{L}'\ots [v_{M}]\ots v_{N}=0,\nn\eeqa
since $\Ome_{\ss{\Mp}}([\ome])=[\ome]$,
$\Ome_{\ss{\Mp}}([J])=[J]$ and $\Ome_{M}([v_{M}])=-[v_{M}]$. Moreover
the Verma module of Virasoro algebra with central charge
$1$ and lowest weight $1$ has a singular vector $2
L(-3)w-4L(-2)L(-1) w+L(-1)^{3}w$ of weight $4$, where $w$ is the cyclic vector
of the Verma module. The image of this singular vector in $\Mm$ is zero, then 
by (\ref{6.4}) and (\ref{6.5}), we have $ g(a_{N},a_{L}) v_{L}'\ots [v_{M}]\ots
v_{N}=0$, where $ g(x,y)=(x-y)(x^{2}-2 x y+ y^{2}-2 x-2
y+1)/2.$

Now suppose that the fusion rule $N_{MN}^{L}$ is nonzero.
Then $v_{L}'\ots [v_{M}]\ots v_{N}$ is nonzero by Proposition \ref{P5.6}
(1). Hence $f(a_{L},a_{N},b_{L})$, $f(a_{N},a_{L},b_{N})$ and
$g(a_{N},a_{L})$ are necessarily zero. By (1) and
Proposition \ref{P4.12}, we assume that $N,L\in\{\Mm,\Ml\
(\lm\neq0),\ \Mt^{\pm}\}$. 

In the cases $N=L=\Mm,\Mtp$ and $\Mtm$ and the cases $L=\Mm$ and
$N=\Mt^{\pm}$, we see that 
$f(a_{L},a_{N},b_{L})$ are nonzero. Therefore the fusion rules of corresponding types
are zero. 

In the case $L=\Mm$ and $N=\Ml$, we see that either $f(1,\lm^{2}/2,-6)$ or
$g(1,\lm^{2}/2)$ is nonzero. Then the fusion rule of corresponding type
is zero. 

 In the case $L=\Ml$ and
$N=\Mu$, assume that the fusion rule is nonzero. Then if we put
$s=\lm^{2}$ and $\mu^{2}=t$, we have 
\eqa &&f\left({s\over 2},{t\over 2},s^{2}-{s\over 2}\right)+ f\left({t\over
2},{s\over 2},t^{2}-{t\over 2}\right)={9\over 16}(s-t)^{2}(3s+3t-2)=0,\nn\\
&&f\left({s\over 2},{t\over 2},s^{2}-{s\over 2}\right)+ f\left({t\over 2},{s\over
2},t^{2}-{t\over 2}\right)+81\,g\left({s\over2},{t\over 2}\right)={9\over
2}(s-t)(s+t-1)=0.\nn
\eeqa
Hence we have $s=t$, that is, $\lm^{2}=\mu^{2}$.
Thus $N_{\ss{ M\Ml}}^{\ss{ \Mu}}=\delta_{\lm^{2},\mu^{2}}$ holds.

In the case $L=\Ml$ and $N=\Mtp$, if put $s=\lm^{2}$, then we have 
\eqa &&f\left({s\over 2},{1\over 16},s^{2}-{s\over 2}\right)={27\over
4096}(8s-1)(32s^{2}-236s+205),\nn\\ 
&&f\left({1\over 16},{s\over
2},{3\over 128}\right)={27\over 4096}(8s-9)(384s^{2}-1160s+131).\nn
\eeqa
But there is no common zero of these two polynomials, and the corresponding fusion
rule is zero.

If $L=\Ml$ and $N=\Mtm$, then we have 
\eqa &&f\left({s\over 2},{9\over 16},s^{2}-{s\over 2}\right)={-9\over
4096}(8s-9)(96s^{2}-996s+119),\nn\\ 
&&f\left({9\over 16},{s\over
2},{-45\over 128}\right)={9\over 8192}(8s-1)(384s^{2}-2504s+2211).\nn
\eeqa
Thus there is no common zero of these two polynomials, and the corresponding fusion
rule is zero in this case too. Now (2) follows from Proposition \ref{P6.2} and
\ref{P6.4}.

{\it Step 3.}  Let $M=\Ml$ ($\lm\neq0$). We prove (3) by dividing into following four
cases: (i) $\lm^{2}\neq 1/2,2,9/2$ case, (ii) $\lm^{2}=2$ case, (iii) $\lm^{2}=9/2$ case,
(iv) $\lm^{2}=1/2$ case.

(i)\ First we assume that $\lm^{2}\neq1/2,$ $2$, $9/2$. Then we can see that
$h(-3)h(-1)\1*v_{M}$ belongs to $L(1,\lm^{2}/2)$ by the direct sum decomposition
(\ref{4.20}). Thus $h(-3)h(-1)\1*v_{M}$ can be expressed by a linear combination of
$L(-m_{1})\cdots L(-m_{k})v_{M}$ ($m_{i}\in\Z_{>0}$,
$m_{1}+\cdots +m_{k}\leq 4$). Then by using (\ref{6.4}) and (\ref{6.5}), we have
following equalities in the contraction $L_{0}^{*}\cdotp A(M)\cdotp N_{0}$:
\eqa &&v_{L}'\ots [(J-4 \ome^{*2}-17\ome +9h(-3)h(-1)\1)*v_{M}]\ots
v_{N}\nn\\
&&\ \ \ =f(a_{L},a_{N},b_{L})v_{L}'\ots [v_{M}]\ots v_{N}=0,\label{6.10}\\
&&v_{L}'\ots [\Ome_{M}(v_{M})*\Ome_{\ss{\Mp}}(J-4 \ome^{*2}-17\ome
+9h(-3)h(-1)\1)]\ots v_{N}\nn\\
&&\ \ \ =v_{L}'\ots \Ome_{M}([(J-4 \ome^{*2}-17\ome
+9h(-3)h(-1)\1)*v_{M}])\ots v_{N}\nn\\
&&\ \ \ =e^{\lm^{2} \pi i\over 2}f(a_{N},a_{L},b_{N})v_{L}'\ots [v_{M}]\ots
v_{N}=0,\label{6.11}
\eeqa
where $f=f(x,y,z)\in\C[x,y,z]$ is given by
\eqa
&&f=z-4 x^{2}+x+{9(\lm^{4}-4
\lm^{2}(x+y)+4(x-y)^{2})\over 32\lm^{2}(\lm^{2}-2)(2\lm^{2}-9)(2\lm^{2}-1)}
\nn\\
&&\ \ \ts(-3(\lm^{2}-2)^{3}+4(8 \lm^{4}-29\lm^{2}+6)
x+4(7\lm^{2}+6)y-4(16\lm^{2}+3)(x-y)^{2}).\nn
\eeqa

We consider the case $N=\Mtp$ and $L=\Mu$ ($\mu\neq0$). If we assume that
the fusion rule $N_{MN}^{L}$ is nonzero, then $v_{L}'\ots [v_{M}]\ots v_{N}$ is
nonzero by Proposition \ref{P5.6} (1). Hence we have
$f(a_{L},a_{N},b_{M})=f(a_{N},a_{L},b_{N})=0$. This shows the following
equations respectively:
\eqa
&&p(t,u)=(8u-1)(8u-9)\nn\\
&&\ \ \ts((1024t+192)u^{2}-(2048t^{2}+512
t+624)u+1024t^{3}-256 t^{2}+816 t+75)=0,\nn\\
&&q(t,u)=((8t+8u-1)^{2}-256tu)\nn\\
&&\ \ \ts((1024t+192)u^{2}-(1024t^{2}-3456
t+816)u+1192 t^{2}+864 t+675)=0.\nn
\eeqa
If $t\in\{1/2,2,9/2\}$, then there is no common solution of
$p(t,u)=q(t,u)=0$. Hence we may assume that $t\notin\{1/2,2,9/2\}$ and discuss by 
interchanging $\lm$ and $\mu$. Thus we have $p(u,t)=q(u,t)=0$. If
$t=1/8,9/8$, then $p(t,u)$ and $q(t,u)$ have no common solution. Therefore $t\neq
1/8$, $9/8$. Similarly $u\neq1/8$, $9/8$. If we put
$p(t,u)=(u-1/8)(u-9/8)r(t,u)$, then $r(t,u)=r(u,t)=0$ hold. Since $r(t,t)$ and
$q(t,t)$ have no common solution, we may also assume that $t\neq u$. From
$r(t,u)-r(u,t)=q(t,u)-q(u,t)=q(t,u)+q(u,t)=0$, we have following three equations:
\eqa &&32\al^{2}-14 \al +45-128\beta=0,\label{6.12}\\
&&(128\beta+3)(64\al^{2}-16\al +1-256\beta)=0,\nn\\
&&(64\al^{2}-16\al +1-256\beta)(64\al^{2}-280\al
+225+1024\beta)=0,\nn 
\eeqa
where we put $\al=t+u$ and $\beta=tu$. If $ 64\al^{2}-16\al
+1-256\beta\neq0$, then $\beta=-3/128$. But In this case, (\ref{6.12}) and
$64\al^{2}-280\al +225+1024\beta=0$ don't have common solution. So we have 
 $64\al^{2}-16\al +1-256\beta=0$ and $32\al^{2}-14 \al +45-128\beta=0$. This
implies $\al=89/12$ and $\beta=30625/2304$. But the solutions
$(t,u)\in\C^{2}$ of the equation $x^{2}-\al x+\beta=0$ don't satisfy $p(t,u)=0$.
Therefore we see that the corresponding fusion rule is zero. 

In the case $N=\Mtm$ and $L=\Mu$
($\mu\neq0$), we can show that the fusion rule $N_{MN}^{L}$ is zero by the
same method in the preceding case. 

We next prove
that $N_{\ss{\Ml\Mu}}^{\ss{\Mn}}=\delta_{\nu^{2},(\lm\pm\mu)^{2}}$ if
$\lm^{2}\neq 1/2$, $2$, $9/2$. By Proposition \ref{P6.1} and Corollary
\ref{C5.8}, it is enough to prove that if $N_{\ss{\Ml\Mu}}^{\ss{\Mn}}\neq 0$ then
$\nu^{2}=(\lm\pm\mu)^{2}$. Let $N=\Mu$ and $L=\Mn$, and assume that
$N_{MN}^{L}$ is nonzero. By Proposition \ref{P5.6}, we see $v_{L}'\ots
[v_{M}]\ots v_{N}\neq 0$. Then the equalities (\ref{6.10}) and (\ref{6.11}) follows
two equations:
\eqa &&(s^{2}+t^{2}+u^{2}-2 st-2 su-2 tu)p(s,t,u)\nn\\
&&{}=(s^{2}+t^{2}+u^{2}-2 st-2
su-2 tu)p(s,u,t)=0,\label{6.17}\eeqa
where $s=\lm^{2}$, $t=\nu^{2}$, $u=\mu^{2}$ and $p\in\C[x,y,z]$ is
given by
\eqa &&p=(-3 x+16xy+z+32yz-z^{2})\nn\\
&&\ \ \ (-2-12x+58xy+16xy^{2}-12 y+3
y^{2}-6 yz-12z+3z^{2}).\nn
\eeqa 
Suppose that $s^{2}+t^{2}+u^{2}-2 st-2 su-2 tu\neq0$, then we have
$p(s,t,u)=p(s,u,t)=0$. In addition, we assume that $\lm^{2}\neq8$. Then the
circle relations $h(-3)h(-1)\1\circ v_{M}$ and $h(-2)^{2}\1\circ v_{M}$
belong to $L(1,\lm^{2}/2)$ and can be expressed by linear combinations of
$L(-m_{1})\cdots L(-m_{k})v_{M}$ ($m_{i}\in\Z_{>0}$, $m_{i}+\cdots
+m_{k}\leq 5$). By (\ref{6.4}) and (\ref{6.5}), we have following equalities:
\eqa &&v_{L}'\ots [h(-3)h(-1)\1\circ v_{M}]\ots
v_{N}=g_{1}(a_{L},a_{N})v_{L}'\ots [v_{M}]\ots v_{N}=0,\label{6.15}\\
&&v_{L}'\ots[h(-2)^{2}\1\circ v_{M}]\ots v_{N}=g_{2}(a_{L},a_{N})v_{L}'\ots
[v_{M}]\ots v_{N}=0,\label{6.16}
\eeqa
where $g_{1}(a_{L},a_{N})$ and $g_{2}(a_{L},a_{N})$ become as follows:
\eqa
&&g_{1}(a_{L},a_{N})=(s^{2}+t^{2}+u^{2}-2 st-2 su-2 tu)(t-u)q(s,t,u),
\nn\\
&&g_{2}(a_{L},a_{N})=(s^{2}+t^{2}+u^{2}-2 st-2 su-2 tu)r(s,t,u),\nn
\eeqa
and $q$, $r\in\C[x,y,z]$ are given by
\eqa
q&=&-12+24x-5x^{2}+12y-16xy+4x^{2}y+3 y^{2}-4 x
y^{2}+12z\nn\\
&&{}-16xz+4x^{2}z+6yz+8xyz-3z^{2}-4yz^{2},\nn\\
r&=&192-245x+108x^{2}-18x^{3}+x^{4}-240y-28xy+8x^{2}y+x^{3}y\nn\\
&&{}+96y^{2}+86xy^{2}-21x^{2}y^{2}-12y^{3}+19xy^{3}-144z+460xz\nn\\
&&{}-152x^{2}z +11x^{3}z-96yz-100xyz+14x^{2}yz+36y^{2}z\nn\\
&&{}-152x
y^{2}z+14xz^{2}+7x^{2}z^{2}+57xyz^{2}-36y^{2}z^{2}+12z^{3}-19xz^{3}.\nn
\eeqa
Since $s^{2}+t^{2}+u^{2}-2 st-2 su-2 tu$ is nonzero, we have $(t-u)q(s,t,u)=r(s,t,u)=0$.
Similarly we have $(u-t)q(s,u,t)=r(s,u,t)=0$ by interchanging $\mu^{2}$ and
$\nu^{2}$. Assume that $t=u$. Then $r(s,t,t)=0$ follows $s=6(1-2t)$, and
$p(6(1-2t),t,t)=0$ implies that $t\notin\{1/2,2,9/2,8\}$. Hence we can interchange
$s$ and $t$, and we have $r(t,6(1-2t),t)=0$ and $r(t,t,6(1-2t))=0$. But the
common solutions of these equations don't satisfy $p(6(1-2t),t,t)=0$. Hence we
see that $t\neq u$, which shows $q(s,t,u)=q(s,u,t)=0$. Next assume that $t=n^{2}/2$
($n=1,2,3,4$). Then $p(s,n^{2}/2,u)=0$ follows 
$$ s={16u^{2}-(16n^{2}+1)u\over 8n^{2}-3}\ \ {\rm or}\ \
-{12u^{2}-12(n^{2}+4)u+3n^{4}-24n^{2}-8\over 4(4n^{2}+29n^{2}-12)}.$$
But in any case, there is no common solution of equations $p(s,u,n^{2}/2)=0$
and $q(s,u,n^{2}/2)=0$. Thus $t\notin\{1/2,2,9/2,8\}$. Similarly we have
$u\notin\{1/2,2,9/2,8\}$. Therefore if we put $x_{1}=s$, $x_{2}=t$ and $x_{3}=u$, 
\eqa
p(x_{i_{1}},x_{i_{2}},x_{i_{3}})=q(x_{i_{1}},x_{i_{2}},x_{i_{3}})
=r(x_{i_{1}},x_{i_{2}},x_{i_{3}})=0\ {\rm and}\ x_{i}\neq x_{j}\ \ {\rm if}\ i\neq
j\label{6.29}\eeqa
 hold for every permutation $\{i_{1},i_{2},i_{3}\}$ of $\{1,2,3\}$ and
$i,j=1,2,3$.  In particular, we have $(t-u)(q(s,t,u)-q(t,s,u))-(t-s)(q(u,t,s)-q(t,u,s))=0$
and this follows
$s+t+u+5=0$ since $s,t,u$ are distinct each other. On the other hand, if we put 
$r(s,t,u)-r(t,s,u)=(s-t)\al(s,t,u)$, then we have $\al(s,t,u)=\al(u,t,s)=0$. If
we next put $\al(s,t,u)-\al(u,t,s)=(s-u)\beta(s,t,u)$, we have
$\beta(s,t,u)=\beta(s,u,t)=0$. Now we see that
$\beta(s,t,u)-\beta(s,u,t)=-16(t-u)(3s+3t+3u-10)=0$ and this follows
$3s+3t+3u-10=0$. It is inconsistent with $s+t+u+5=0$. Thus we see that there is no
solution satisfy (\ref{6.29}). Hence $s^{2}+t^{2}+u^{2}-2 st-2 su-2 tu=0$. By
substituting
$\lm^{2}$ for $s$, $\nu^{2}$ for $t$ and $\mu^{2}$ for $u$, we have
$\nu^{2}=(\lm\pm\mu)^{2}$ if
$\lm\neq1/2$, $2$, $9/2$, $8$. 

In the case $\lm^{2}=8$, we may assume that
$N$ and $L$ are any of $\Mu$ ($\mu^{2}=1/2,2,9/2,8$), that is to say, that $\mu^{2}$
and $\nu^{2}$ are any of $1/2$, $2$, $9/2$ and $8$. But then
$p(8,\mu^{2},\nu^{2})$ is nonzero. Hence we have $\nu^{2}=(\mu\pm2\sqrt{2})^{2}$
by (\ref{6.17}). Consequently we see that (3) hold if $\lm^{2}\neq1/2$, $2$,
$9/2$.

(ii) We next consider the case $\lm^{2}=2$. Then we have the irreducible
decomposition for Virasoro algebra $\Ml=L(1,1)\ops L(1,4)\ops L(1,9)\ops\cdots$
(see (\ref{4.20})). Put
$u=\sqrt{2}h(-3)v_{M}-3h(-2)h(-1)v_{M}+\sqrt{2}h(-1)^{3}v_{M}$ which is the
lowest weight vector of $L(1,4)$. Since
$h(-3)h(-1)\1*v_{M}\in L(1,1)\ops L(1,4)$, it is a linear
combination of $L(-1)u$, $u$ and $L(-m_{1})\cdots L(-m_{k})v_{M}$
$(m_{i}\in\Z_{>0}$, $m_{1}+\cdots+m_{k}\leq 4$). Then by formulas (\ref{6.4}) and
(\ref{6.5}), we have following equation:
\eqa &&v_{L}'\ots [(J-4 \ome^{*2}-17\ome +9h(-3)h(-1)\1)*v_{M}]\ots
v_{N}\nn\\
&&\ \ \ =f_{1}(a_{L},a_{N})v_{L}'\ots [u]\ots
v_{N}+f_{2}(a_{L},a_{N},b_{L})v_{L}'\ots [v_{M}]\ots v_{N}=0,\nn
\eeqa
where $f_{1}=f(x,y)\in\C[x,y]$ and $f_{2}=f(x,y,z)\in\C[x,y,z]$ are given by
$$
f_{1}={3\over 2}+{21\over 8}(x-y),$$
$$f_{2}=z+{95\over 8}x-{99\over 8}y-{173\over
8}x^{2}-{9\over4}xy+{207\over 8}y^{2}+{47\over
4}x^{3}-27x^{2}y+{135\over 4}xy^{2}-{27\over 2}y^{3}.$$
Since $\Ome_{M}([u])=[u]$ and
$\Ome_{M}([v_{M}])=-[v_{M}]$, we have 
$$
-f_{1}(a_{N},a_{M})v_{N}'\ots [u]\ots v_{N}+f_{2}(a_{N},a_{L},b_{N})v_{L}'\ots
[v_{M}]\ots v_{N}=0.$$

We may assume that $N$ and $L$ are any of $\Mu$ ($\mu^{2}=1/2,2,9/2$),
$\Mt^{\pm}$. But then we can see that the determinant of the matrix
\eqa\left(\begin{array}{cc}
f_{1}(a_{L},a_{N})&f_{2}(a_{L},a_{N},a_{L}) \cr
-f_{1}(a_{N},a_{L})&f_{2}(a_{N},a_{L},a_{N})
\end{array}\right)\nn\eeqa
is nonzero except for the following cases:
\eqa (N,L)&=&(\Mu,\Mn);\ \
(\mu^{2},\nu^{2})=(1/2,9/2),(9/2,1/2),(1/2,1/2),\nn\\
&&{}(\Mt^{\al},\Mt^{\beta});\,\,\al,\beta\in\{+,-\}.\nn
\eeqa
Note that $9/2=(\sqrt{2}+1/\sqrt{2})^{2}$ and
$1/2=(\sqrt{2}-1/\sqrt{2})^{2}=(\sqrt{2}-3/\sqrt{2})^{2}$. Then we have 
$ N_{\ss{\Ml\Mu}}^{\ss{\Mn}}=\delta_{\nu^{2},(\lm\pm\mu)^{2}},\ \
N_{\ss{\Ml\Mu}}^{\ss{\Mt^{\pm}}}=0$
for every $\mu^{2},\nu^{2}=1/2,2,9/2$ by Corollary \ref{C5.8} and Proposition
\ref{P6.1}. Together with the results of (i) and Proposition
\ref{P6.4}, we see that (3) holds for $\lm^{2}=2$.

(iii)  We consider the case $\lm^{2}=9/2$. In this case, the submodule
for Virasoro algebra of $M=\Ml$ generated by $v_{M}$ is isomorphic to
the irreducible module $L(1,9/4)$. The Verma module for Virasoro algebra with central
charge $1$ and lowest weight $9/4$ has a singular vector $18L(-4)w-14L(-3)L(-1)
w-9L(-2)^{2}w+10L(-2)L(-1)^{2}w-L(-1)^{4}w$, where $w$ is the cyclic
vector of the Verma module. Since the image of the singular vector in
$\Ml$ is zero,  by using
(\ref{6.4}) and (\ref{6.5}), we have the following equality:
$$ f(a_{L},a_{N}) v_{L}'\ots [v_{M}]\ots v_{N}=0,$$
where $f=f(x,y)\in\C[x,y]$ is given by
$$f=(81-72(x+y)+16(x-y)^{2})(1-8(x+y)+16(x-y)^{2}).$$
By the results of Step 1, Step 2, and (i), (ii) of Step 3, we may assume that $N$
and $L$ are any of $\Mu$ ($\mu^{2}=1/2,9/2$), $\Mt^{\pm}$. If
$a_{L},a_{N}\in\{1/4,9/4,1/16,9/16\}$, then $ f(a_{L},a_{N})$ is nonzero except
for the pairs 
$(a_{L},a_{N})=$$(1/16,1/16)$, $(1/16,9/16)$, $(9/16,1/16)$, $(9/16,9/16)$. By
Corollary
\ref{C5.8} and Proposition \ref{P6.1}, we see that for
$\lm^{2}=9/2$, $\mu^{2},\nu^{2}=1/2,9/2$, $
N_{\ss{\Ml\Mu}}^{\ss{\Mn}}=N_{\ss{\Ml\Mu}}^{\ss{\Mt^{\pm}}}=0$. Then
Proposition \ref{P6.4} implies that (3) holds for $\lm^{2}=9/2$.

(iv) We prove that (3) holds for $\lm^{2}=1/2$. In this case we have the irreducible
decomposition $\Ml=L\left(1,{1/4}\right)\ops L\left(1,{9/4}\right)\ops
L\left(1,{25/4}\right)\ops\cdots$  (see (\ref{4.20})). 
Let $u=\sqrt{2}h(-2)v_{M}-2h(-1)^{2}v_{M}$ which is the lowest weight
vector of $L(1,9/4)$. Then $h(-3)h(-1)\1*v_{M}$ belongs to $L(1,1/4)\ops L(1,9/4)$
and is a linear combination of $L(-m_{1})\cdots L(-m_{k})v_{M}$ and
$L(-n_{1})\cdots L(-n_{\ell})u$ ($m_{i},n_{j}\in\Z_{>0}$, $m_{1}+\cdots+m_{k}\leq
4,$ $n_{1}+\cdots+n_{\ell}\leq 3$). By (\ref{6.4}) and (\ref{6.5}), we have the
following equality:
\eqa &&v_{L}'\ots [(J-4 \ome^{*2}-17\ome +9h(-3)h(-1)\1)*v_{M}]\ots
v_{N}\nn\\
&&\ \ \ =f_{1}(a_{L},a_{N})v_{L}'\ots [u]\ots
v_{N}+f_{2}(a_{L},a_{N},b_{L})v_{L}'\ots [v_{M}]\ots v_{N}=0,\label{6.18}
\eeqa
where $f_{1}=f(x,y)\in\C[x,y]$ and $f_{2}=f(x,y,z)\in\C[x,y,z]$ are given by
\eqa
f_{1}={27\over 128}+{13\over 16}x-{19\over 16}y+{11\over
16}(x-y)^{2},\ \ f_{2}=z-{3\over 2}+12(x+y)-24(x-y)^{2}.\nn\eeqa
Since $\Ome_{M}([u])=e^{\pi i/4}[u]$ and $\Ome_{M}([v_{M}])=e^{\pi i/4}[v_{M}]$, we
have also
\eqa
f_{1}(a_{N},a_{L})v_{L}'\ots [u]\ots
v_{N}+f_{2}(a_{N},a_{L},b_{N})v_{L}'\ots [v_{M}]\ots v_{N}=0.\nn
\eeqa
In addition, we have $L(-1)^{2}v_{M}+L(-2)v_{M}=0$, since the Verma module
of central charge $1$ and highest weight $1/4$ has nontrivial singular vector of
weight $9/4$. Therefore we have the following equality:
\eqa v_{L}'\ots [(L(-1)^{2}v_{M}+L(-2)v_{M})]\ots
v_{N}=g(a_{L},a_{N})v_{L}'\ots [v_{M}]\ots v_{N}=0,\label{6.20}
\eeqa
where $g=g(x,y)=-1/16+(x+y)/2-(x-y)^{2}$. By Step 1, Step 2 and (i)-(iii) of Step 3, we
may assume that $N$ and $L$ are any of $\Ml$ ($\lm^{2}=1/2$), $\Mt^{\pm}$. 

In the case $N=L=\Ml$ ($\lm^{2}=1/2$), since both $f_{1}(a_{L},a_{N})$
and $g(a_{L},a_{N})$ are nonzero, so $v_{L}'\ots [u]\ots v_{N}=v_{L}'\ots
[v_{M}]\ots v_{N}=0$ by (\ref{6.18}) and (\ref{6.20}). Hence we have 
$\dim L_{0}^{*}\cdotp A(M)\cdotp N_{0}=0$ by Proposition \ref{P5.6}. Thus
Proposition \ref{P3.12} shows that the fusion rule $N_{MN}^{L}$ is zero.

In the case $L=\Mu$ ($\mu^{2}=1/2$) and $N=\Mt^{\pm}$, the determinant of
the matrix
$$\left(\begin{array}{cc}
f_{1}(a_{L},a_{N})&f_{2}(a_{L},a_{N},a_{L}) \cr
f_{1}(a_{N},a_{L})&f_{2}(a_{N},a_{L},a_{N})
\end{array}\right)$$
is nonzero. Hence by Proposition \ref{P3.12}, the fusion rule $N_{MN}^{L}$ is zero in
this case too.

In the case that $N$ and $L$ are either $\Mtp$ or $\Mtm$, 
$v_{L}'\ots [u]\ots v_{N}$ and $v_{L}'\ots[v_{M}]\ots v_{N}$ are linearly
dependent, so the dimension of $L_{0}^{*}\ots A(M)\ots N_{0}$ is less than one by
Proposition \ref{P5.6} (3).  Now Proposition \ref{P6.4} (1) shows that the fusion rule
$N_{MN}^{L}$ is 1 in this case. Consequently, (3) holds for $\lm^{2}=1/2$. Thus we see
that (3) holds for all $\lm\in\C-\{0\}$.

{\it Step 4.}  In the case $M=\Mtp$, we
have $\Mtp = L\left(1,{1/16}\right) \ops L\left(1,{49/16}\right) \ops 
L\left(1,{81/16}\right) \ops \cdots$  (see (\ref{4.31})). We put
$u=9h(-5/2)h(-1/2)1-5 h(-3/2)^{2}v_{M}-10h(-3/2)h(-1/2)^{3}+4h(-1/2)^{6}1$
which is the lowest weight vector of $L(1,49/16)$. Since $h(-3)h(-1)\1*v_{M}$ is in
$L(1,1/16)\ops L(1,49/16)$, it can be expressed a linear combination of $L(-1)u$,
$u$ and 
$L(-m_{1})\cdots L(-m_{k})v_{M}$ ($m_{i}\in\Z_{>0}$,
$m_{1}+\cdots+m_{k}\leq 4$). So we have following equality:
\eqa &&v_{L}'\ots [(J-4 \ome^{*2}-17\ome +9h(-3)h(-1)\1)*v_{M}]\ots
v_{N}\nn\\
&&\ \ \ =f(a_{L},a_{N})v_{L}'\ots [u]\ots
v_{N}+g(a_{L},a_{N},b_{L})v_{L}'\ots [v_{M}]\ots v_{N}=0,\nn
\eeqa
where $f=f(x,y)\in\C[x,y]$ and $g=g(x,y)\in\C[x,y,z]$ are given by
\eqa f={1\over 2}+{8\over 7}(x-y),&&g=5z-{135\over 1792}-{1\over
56}x+{73\over 28}y-{82\over 7}x^{2}+{212\over 7}xy-{180\over 7}
y^{2}\nn\\
&&\ \ \ \ \ +{32\over 7}(x-y)^{2}(5x+12y)-{256\over
7}(x-y)^{4}.\nn\eeqa Since $\Ome_{M}([u])=-e^{\pi i/16}[u]$ and
$\Ome_{M}([v_{M}])=e^{\pi i/16}[v_{M}]$, we have equality 
$$
-f(a_{L},a_{N})v_{L}'\ots [u]\ots
v_{N}+g(a_{L},a_{N},b_{L})v_{L}'\ots [v_{M}]\ots v_{N}=0.$$

By the results of Step 1-Step 3 and Proposition \ref{P4.12}, we assume
that $N$ and $L$ are either $\Mtp$ or $\Mtm$. But then the determinant of the
matrix 
\eqa\left(\begin{array}{cc}
f(a_{L},a_{N})&g(a_{L},a_{N},a_{L}) \cr
-f(a_{N},a_{L})&g(a_{N},a_{L},a_{N})
\end{array}\right)\nn\eeqa
is nonzero. Hence we have $v_{L}'\ots [u]\ots v_{N}=v_{L}'\ots [v_{M}]\ots
v_{N}=0$. This proves that $N_{\ss{\Mtp\Mt^{\al}}}^{\ss{\Mt^{\beta}}}=0$ for any
$\al,\beta\in\{+,-\}$. Thus we see that (4) holds.

{\it Step 5.}   By the results of Step 1-Step 4, it is enough to show that
$N_{MN}^{L}=0$ for $M=N=L=\Mtp$. Since we have the direct product decomposition 
$\Mtm=L\left(1,{9/16}\right)\ops L\left(1,{25/16}\right)\ops
L\left(1,{121/16}\right)\ops\cdots$ (see \ref{4.32}),
if we put $u=-(1/2)h(-3/2)1+h(-1/2)^{3}1$ which is the lowest
weight vector of $L(1,25/16)$, then $h(-3)h(-1)\1*v_{M}$  can be expressed a linear
combination of 
$L(-m_{1})\cdots L(-m_{k})v_{M}$ and $L(-n_{1})\cdots L(-n_{\ell})u$
($m_{i},n_{j}\in\Z_{>0}$, $m_{1}+\cdots+m_{k}\leq 4,$ $n_{1}+\cdots+n_{\ell}\leq
3$). Calculating the vector $v_{L}'\ots [(J-4 \ome^{*2}-17\ome
+9h(-3)h(-1)\1)*v_{M}]\ots v_{N}$ by means of (\ref{6.4}) and (\ref{6.5}), we have
the following linear equalities: 
\eqa {75\over 224}v_{L}'\ots [u]\ots v_{N}-{135\over 256}v_{L}'\ots
[v_{M}]\ots v_{N}=0.\nn\eeqa
Since $\Ome_{M}([u])=-e^{9\pi i/16}[u]$ and $\Ome_{M}([v_{M}])=e^{9\pi
i/16}[v_{M}]$, by Proposition \ref{P3.5}, we have
\eqa -{75\over 224}v_{L}'\ots [u]\ots v_{N}-{135\over 256}v_{L}'\ots
[v_{M}]\ots v_{N}=0.\nn\eeqa
This follows that $ v_{L}'\ots [u]\ots v_{N}=v_{L}'\ots [v_{M}]\ots
v_{N}=0$. So we see that the fusion rule $N_{MN}^{L}=0$ if
$M=N=L=\Mtp$ by Proposition \ref{P5.6} (3). The proof of the theorem is
complete.\qed
\vs

{\bf Aknoledgments :} I would like to thank Professor Kiyokazu Nagatomo for useful
discussions and suggestions. I also thank Doctor Yoshiyuki Koga and Akihiko Ogawa for
many opinions.

\end{document}